\documentclass[12pt,letterpaper,reqno]{amsart}

\usepackage{epsfig}
\usepackage{amsmath}
\usepackage{amssymb}
\usepackage{amsthm}
\usepackage{indentfirst}
\usepackage{xspace}
\usepackage{multirow}
\usepackage{hyperref}
\usepackage{xcolor}
\definecolor{darkred}{rgb}{0.2,0.25,0.75}
\hypersetup{colorlinks=true,urlcolor=darkred,linkcolor=darkred,citecolor=darkred}
\usepackage{verbatim}
\usepackage[letterpaper,margin=1in]{geometry}
\usepackage{mathpazo}
\usepackage{thmtools}
\usepackage[all]{xy}
\usepackage{longtable}
\usepackage{capt-of}
\usepackage{enumitem}

\setlist{itemsep = 0.20em, topsep = 0.20em}

\declaretheoremstyle[spaceabove=0.25cm,spacebelow=0.25cm,notefont=\normalfont\bfseries, notebraces={(}{)}]{theorem}
\declaretheoremstyle[spaceabove=0.25cm,spacebelow=0.25cm,bodyfont=\normalfont,notefont=\normalfont\bfseries, notebraces={(}{)}]{noital}
\declaretheoremstyle[spaceabove=0.25cm,spacebelow=0.25cm,bodyfont=\normalfont\color{darkgreen},notefont=\normalfont\bfseries, notebraces={(}{)}]{green}
\declaretheoremstyle[spaceabove=0.25cm,spacebelow=0.25cm,bodyfont=\normalfont,notefont=\normalfont\bfseries,qed=$\qedsymbol$,notebraces={(}{)}]{proofstyle}

\declaretheorem[name=Theorem,style=theorem]{thm}
\declaretheorem[name=Question,style=theorem]{question}
\declaretheorem[name=Conjecture,style=theorem]{conj}

\numberwithin{equation}{section}

\newcommand{\fA}{{\mathfrak A}}

\newcommand{\cC}{\ensuremath{\mathcal C}}

\newcommand{\cM}{\ensuremath{\mathcal M}}

\newcommand{\cP}{\ensuremath{\mathcal P}}

\newcommand{\cO}{\ensuremath{\mathcal O}}

\newcommand{\cX}{\ensuremath{\mathcal X}}

\newcommand{\cW}{\ensuremath{\mathcal W}}

\newcommand{\R}{\ensuremath{\mathbb R}}
\newcommand{\C}{\ensuremath{\mathbb C}}
\newcommand{\PP}{\ensuremath{\mathbb P}}
\newcommand{\Z}{\ensuremath{\mathbb Z}}

\newcommand{\half}{\ensuremath{\frac{1}{2}}}

\newcommand{\hk}{hyperk\"ahler\xspace}
\newcommand{\Hk}{Hyperk\"ahler\xspace}

\newcommand{\I}{{\mathrm i}}
\newcommand{\e}{{\mathrm e}}
\newcommand{\de}{\mathrm{d}}

\newcommand{\here}{{\mathrm {here}}}
\newcommand{\there}{{\mathrm {there}}}
\newcommand{\SL}{{\mathrm {SL}}}
\newcommand{\SO}{{\mathrm {SO}}}
\newcommand{\Gr}{{\mathrm {Gr}}}

\newcommand{\abs}[1]{\lvert#1\rvert}

\newcommand{\IP}[1]{\langle#1\rangle}

\newcommand{\ti}[1]{\textit{#1}}

\DeclareMathOperator{\re}{Re}
\DeclareMathOperator{\Tr}{Tr}

\DeclareMathOperator{\Hom}{Hom}
\DeclareMathOperator{\End}{End}

\DeclareMathOperator{\rank}{rank}
\DeclareMathOperator{\genus}{genus}
\DeclareMathOperator{\holes}{holes}

\newcommand{\insfig}[3]{\begin{figure}[htbp] \centering \includegraphics[scale=#2]{figures/#1.pdf} \caption{#3} \label{fig:#1} \end{figure}}

\begin{document}

\bibliographystyle{utphys}

\setcounter{page}{1}

\title{Integral iterations for harmonic maps}
\author[A. Neitzke]{Andrew Neitzke}
\date{}

{\abstract{We study minimal harmonic maps $g: \C \to \SO(3) \backslash \SL(3,\R)$, parameterized by polynomial cubic differentials $P$ in the plane. The asymptotic structure of such a $g$ is determined by a convex polygon $Y(P)$
in $\R\PP^2$. We give a
conjectural method for determining $Y(P)$
by solving a fixed-point problem for a certain integral operator. 
The technology of spectral networks and BPS state counts
is a key input to the formulation of this fixed-point problem.
We work out two families of examples in detail.
}
}}

\maketitle


\section{Introduction}

\subsection{Summary}

This paper concerns harmonic maps from the flat 
complex plane to a Riemannian symmetric space,
\begin{equation} \label{eq:g-map-intro}
  g: \C \to \SO(3) \backslash \SL(3,\R),
\end{equation}
which are \ti{minimal} and have \ti{polynomial growth}.
The basic facts about such maps were set out in
\cite{MR3432157}:
as we review in \S\ref{sec:harmonic-maps} below,
\begin{itemize}
\item for each degree $n$ polynomial $P$ there is a corresponding
harmonic map $g$, unique
up to an $\SL(3,\R)$ translation of $\SO(3) \backslash \SL(3,\R)$,
\item the asymptotic
behavior of $g(z)$ as $z \to \infty$ is controlled by
a convex $(n+3)$-gon in $\R\PP^2$, which we call $Y(P)$.
\end{itemize}
Thus the question arises: 

\begin{question} \label{ques:harmonic}
What can we say about
the convex $(n+3)$-gon $Y(P)$ corresponding to a given
degree $n$ polynomial $P$?
\end{question}

The main aim of this paper is to describe a conjectural approach to
this question which arises naturally from
the work
\cite{Gaiotto:2008cd,Gaiotto:2009hg,Gaiotto2012}.\footnote{More precisely, this paper involves a small extension of \cite{Gaiotto2012} to treat connections with a specific sort of irregular singularity.}
For the most part we do \ti{not} attempt to explain the motivation, 
focusing rather on giving enough detail
to make it clear precisely what the conjecture is.
In short:
\begin{conj} \label{conj:main}
When $P$ has only simple zeroes, $Y(P)$ can be computed
from the solution of a fixed-point problem for
an integral operator acting on collections of functions in
one variable $\zeta$, described in \S\ref{sec:integral-equations}
below.
\end{conj}

The fixed-point problem appearing in \autoref{conj:main} 
does not involve the polynomial $P$ directly
but rather uses various ingredients derived from
$P$, which are introduced in \S\ref{sec:background}.
In particular, it uses the \ti{WKB spectral networks}
derived from $P$, described in \S\ref{sec:spectral-network}.

In \S\ref{sec:case-studies} we work out the details of
two examples of \autoref{conj:main}:
one with $n=2$ (in \S\ref{sec:pentagon}-\S\ref{sec:network-analysis-pentagon}), another with $n=3$ (in \S\ref{sec:hexagon}-\S\ref{sec:network-analysis-hexagon}).

\autoref{conj:main} 
is closely related to an integral 
iteration introduced much earlier in the context of $tt^*$ geometry
\cite{MR1213301,Cecotti:1993rm} and to the thermodynamic Bethe ansatz.
It is complementary
to the more obvious method of solving the harmonic
map equations directly, and exposes 
different aspects of the underlying structure.
Its consequences are especially transparent when we study
the asymptotic behavior for large $P$; this is explained
in \S\ref{sec:exact}-\S\ref{sec:asymptotics}.

\subsection{A sample prediction} \label{sec:sample}
Here is a concrete example (from \S\ref{sec:pentagon}) of a prediction obtained from \autoref{conj:main}.
We consider a $1$-parameter family of polynomials:
\begin{equation} \label{eq:pentagon-intro}
  P(z) = \frac{R^3}{2} \left(-z^2 + 1\right), \qquad R \in \R_+.
\end{equation}
The family \eqref{eq:pentagon-intro} corresponds to a $1$-parameter family of convex pentagons 
$Y(P)$ in $\R\PP^2$,
up to $\SL(3,\R)$ action.
Such a pentagon is determined by two
invariant cross-ratios, which we denote $(X_{\gamma_1},X_{\gamma_2})$.\footnote{For the definition of $X_{\gamma_i}$ in terms
of the vertices of the pentagon, see \eqref{eq:pentagon-coords} below.} 
The fact that the coefficients in
\eqref{eq:pentagon-intro} are all real leads to
a reflection symmetry which implies $X_{\gamma_2}(R) = 1$ for all $R$,
so only $X_{\gamma_1}(R)$ remains to be determined.
We predict (see \eqref{eq:asymptotics-specialize-pentagon} below)
\begin{equation} \label{eq:asymptotics-specialize-pentagon-intro}
  X_{\gamma_1}(R) = \exp \left( aR - \frac{3}{2 \sqrt{\pi \rho R}} \e^{-2 \rho R} + \delta(R) \right), \quad   a \approx -4.00648, \quad \rho \approx 2.31315,
\end{equation}
where the remainder function $\delta(R)$ obeys
\begin{equation}
 \lim_{R \to \infty} \delta(R) \sqrt{R} \e^{2 \rho R} = 0.
\end{equation}
Two features of \eqref{eq:asymptotics-specialize-pentagon-intro}
deserve comment:
\begin{itemize}
\item The leading behavior of $\log X_{\gamma_1}(R)$ is \ti{linear} in $R$.
The constant $a$ arises as a period:
$aR = 2 \re \oint_{\gamma_1} P(z)^{1/3} \de z$, where the contour
$\gamma_1$ is shown in Figure \ref{fig:pentagon-cover-cycles}
below.
\item The subleading behavior of $\log X_{\gamma_1}(R)$ is \ti{exponentially suppressed}
in $R$, and the coefficient of the exponentially suppressed term
is explicitly computed.
\end{itemize}
These two features are aspects of the general \autoref{conj:main},
which recur more generally.
In \eqref{eq:asymptotics-specialize} below 
we give the general form of an
asymptotic prediction generalizing \eqref{eq:asymptotics-specialize-pentagon-intro}.

\subsection{Numerical checks} \label{sec:numerical} It is possible to subject
\autoref{conj:main} to direct checks. Indeed, \autoref{conj:main}
leads to a numerical scheme for computing the invariants of $Y(P)$,
by iteration of the relevant integral operator.
On the other hand, one can compute
the invariants of $Y(P)$ directly, by solving the harmonic 
map equations numerically; tools for doing so
have been developed by David Dumas and Michael Wolf \cite{blaschke}.

In joint work in progress,
Dumas and I apply both of these methods
to the examples in \S\ref{sec:case-studies}, and compare
the results.\footnote{Preliminary results are encouraging, but not conclusive. 
One sample (cherry-picked) data point: in the example of \S\ref{sec:sample},
our current estimate of
$X_{\gamma_1}(R = 0.5)$ is $\approx 0.1286$, the two methods of computation
agreeing to this precision.} 
The possibility of making such a comparison was the main motivation
for writing this paper.
Moreover, preliminary such comparisons
were important in shaking out various mistakes
which I made in initial attempts to write \S\ref{sec:case-studies}.

\subsection{Contours and clusters} In the example of \S\ref{sec:sample}
above, the 
asymptotics of the cross-ratio $X_{\gamma_1}$ are controlled by a specific
period integral over a contour $\gamma_1$. This is another feature 
of the general
\autoref{conj:main}: there is a lattice $\Gamma$ of allowed contours, and for
each $\gamma \in \Gamma$, there is a corresponding cross-ratio\footnote{I abuse terminology by using ``cross-ratio'' for any $\SL(3,\R)$-invariant
function of the vertices of $Y(P)$.} $X_\gamma$,
such that the leading asymptotics of $X_\gamma$ are controlled by the period
$\oint_\gamma P(z)^{1/3} \de z$. Moreover, knowing all 
the cross-ratios $X_\gamma$ is sufficient to determine the polygon $Y(P)$.

One of the essential points in making \autoref{conj:main} precise 
is thus to identify which cross-ratio $X_\gamma$
corresponds to a given contour $\gamma$. 
This turns out to be rather subtle: to do it, we need
the spectral network associated to the polynomial $P$.
Given the network we face the problem of drawing certain
compatible \ti{abelianization trees} which determine the $X_\gamma$.
We have solved this problem by hand for the
examples we treat in this paper, but it would be very desirable
to have a more general method, or at least to show that a solution
always exists.
This part of the story is described in \S\ref{sec:spectral-network}-\S\ref{sec:combining-abelianization-trees}.

The underlying combinatorics appears to be closely 
connected to the theory of cluster algebras, specifically
the cluster structure on the homogeneous coordinate ring of $\Gr(3,n+3)$:
roughly, each spectral network corresponds to a cluster.
For a brief account of this see \S\ref{sec:cluster}, but there
is much more to be understood here.

\subsection{Non-simple zeroes}
It would be very interesting to understand the generalization
of \autoref{conj:main} when the zeroes of $P(z)$ are not required
to be simple.
In the case where exactly two zeroes collide, a generalization
of the integral operator from \autoref{conj:main} 
written in \cite{Garz2017} may be relevant.
Indeed, the main example considered there (``pentagon example'')
is directly related to the case $n=2$ of this paper.

\subsection{Other ranks} \label{sec:other-ranks} 
There is a natural generalization of 
\eqref{eq:g-map-intro} to harmonic maps
\begin{equation}
  g: \C \to \SO(K) \backslash \SL(K,\R),
\end{equation}
again with a polynomial growth condition
labeled by an integer $n$. In this case the maps $g$ are labeled
by \ti{tuples} of polynomials $(P_2, P_3, \dots, P_K)$, and the condition
is that the roots of $\lambda^K + \sum_{i=2}^K \lambda^{K-i} P_i(z) = 0$
should behave locally 
as $\lambda_a \sim c_a z^{\frac{n}{K}}$ at large $\abs{z}$, with a 
different constant $c_a \in \C$ for each root.\footnote{Strictly speaking, this is a generalization of our setup even when $K=3$, since it allows two polynomials
$P_2, P_3$ instead of just $P = P_3$, cf. \eqref{eq:our-special-case}. 
This generalization corresponds to taking general harmonic maps instead of minimal ones.}\footnote{Incidentally, when $K \vert n$,
the constants $c_a \in \C$ have a global meaning,
and varying them while holding them distinct
should induce an action of the braid group $B_K$ on
the moduli space of flat irregular connections
by ``iso-Stokes deformation.'' The fact that
such group actions can come from variation of parameters at 
irregular singularities has been emphasized by Philip Boalch;
see e.g. \cite{Boalch} and references therein.
It would be interesting to know whether this is 
the origin of the $B_K$ action on $\Gr(K,rK)$
recently found in \cite{Fraser2017}.
}
The asymptotic behavior of such a map should be
controlled by a point of $\Gr_\R(K, n+K) / (\R^\times)^{n+K-1}$.
\autoref{conj:main} should extend directly to this situation.

In the case $K=2$ these harmonic maps have indeed been studied,
e.g. in \cite{MR1362649}.
Moreover, the extension of \autoref{conj:main} to $K=2$ is actually
simpler than the case $K=3$ which we treat in this paper.
This is because the combinatorics of spectral networks
in this case is easier: each spectral network
gives a triangulation of the $(n+2)$-gon, and one reads off
the cross-ratios $X_\gamma$ directly 
from this triangulation. 
The necessary rules were worked out in \cite{Gaiotto:2009hg}.\footnote{
Using these rules, an application of the $K=2$ version 
of \autoref{conj:main} to a closely related
problem, concerning minimal surfaces in $AdS_3$ with polygonal boundary, was already made in \cite{Alday:2009yn}. See also \cite{Alday:2009dv} where this
was extended to a special case of minimal surfaces in $AdS_5$.}
Thus, it might have made 
sense to write this paper about $K=2$ rather than $K=3$.
The main reason for choosing instead $K=3$ is the
accidental fact that the tools 
\cite{blaschke} were designed specifically for that case.

Working out concrete examples with $K>3$, on the other hand, is likely 
infeasible without some improvement
in our understanding of the combinatorics of spectral networks.

\subsection{Real and complex Higgs bundles} \label{sec:higgs-bundles} One reason
for our interest in \autoref{ques:harmonic} is that it is 
a special case (indeed, almost the simplest nontrivial case) of 
the \ti{nonabelian Hodge correspondence} which relates
Higgs bundles and flat connections:

\begin{question} \label{ques:higgs-real} 
What can we say about the invariants of the flat $\SL(K,\R)$-connection
corresponding to a given rank $K$ real Higgs bundle over a 
Riemann surface $C$, with or without singularities at points of $C$?
\end{question}

\autoref{ques:harmonic} is the special case of \autoref{ques:higgs-real} where we choose $K=3$, $C = \C\PP^1$,
fix a specific sort of irregular singularity at $z = \infty$,
and restrict to Higgs bundles with $\Tr \varphi^2 = 0$.
We explain this in \S\ref{sec:nonabelian-hodge}.
In turn, \autoref{ques:higgs-real} is a special case of a similar
question for general Higgs bundles and complex flat connections.
A version of \autoref{conj:main} is expected to
apply to these more general questions as well; 
this was the original
context of \cite{Gaiotto:2008cd,Gaiotto:2009hg,Gaiotto2012}.\footnote{In 
sufficiently complicated cases ---
particularly the case of Higgs bundles without singularities --- the WKB spectral networks are \ti{dense}
on $C$. This case is more difficult, since problems of \ti{dynamics}
mix with the combinatorial problems. Some of the
necessary tools in case $K=2$ have been developed in
\cite{Fenyes2015}.}

\subsection{\Hk metrics} Moduli spaces of (complex) Higgs bundles with irregular
singularities, appropriately defined, 
are expected to carry natural \hk metrics, much like their
cousins without singularities.
In this paper, we do not explore these
metrics. For context, though, we remark that
the functions $\cX_\gamma(\zeta)$
computed by the integral iteration in \S\ref{sec:integral-equations}
are the main ingredient in a conjectural recipe for constructing these metrics, on the specific moduli spaces of Higgs bundles with irregular
singularity appearing
here; this recipe is described in \cite{Gaiotto:2008cd}. 
The parameter $\zeta$ arises as the coordinate on
the twistor sphere.

\subsection{The WKB method}

The problem of determining the large-$P$ behavior of $Y(P)$ 
is part of a broader class of asymptotic
problems which have been investigated at some length.
I cannot really do this history justice, but at least some
of the key references are: 
\begin{itemize}
\item \cite{MR1133870,MR1152229,MR2350148,Katzarkov2013,Collier2014,Mazzeo2014,Katzarkov2015,Mochizuki2015} for Higgs bundles over a compact Riemann surface,
\item \cite{MR1417948,MR1362649,MR3432157} for harmonic maps $\C \to \SO(K) \backslash \SL(K,\R)$ for $K = 2$ or $K = 3$,
\item \cite{MR2705981,Acost2016} for opers over 
a compact Riemann surface.
\end{itemize}
A recurring theme in this area is the \ti{WKB method}
for studying families of flat connections,
\begin{equation}
  \nabla(t) = t^{-1} \varphi + D + t A_1 + \cdots,
\end{equation}
by studying flat sections in an expansion around $t = 0$.
This is often applied taking $t = 1/R$, where $R$
rescales the Higgs field or the polynomial $P$.
It seems possible that the
leading term in our asymptotic predictions,
e.g. the $aR$ in \eqref{eq:asymptotics-specialize-pentagon-intro},
could be obtained rigorously by a careful application of this technique,
using the spectral network to organize the patching-together of local 
estimates.
The work \cite{Katzarkov2013,Katzarkov2015}, relating
the spectral network to an asymptotic map to a building,
seems to be connected to this picture.

One motivation of \autoref{conj:main}, described 
in \cite{Gaiotto:2009hg,Gaiotto2012},
also involves the WKB method, but 
rather than taking $t = 1/R$ one considers
a family of flat connections $\nabla(\zeta)$, $\zeta \in \C^\times$, 
associated to a \ti{fixed} $P$,
and takes $t = \zeta$. (The relevant family appears as \eqref{eq:connection-family} below.)

\subsection{Quantum field theory} The original motivation of
\cite{Gaiotto:2008cd,Gaiotto:2009hg,Gaiotto2012}
was to address questions about BPS states in supersymmetric
quantum field theories of ``class $S$.''
The specific problem which we study here is related to
the \ti{(generalized) Argyres-Douglas theories} of \ti{type $(A_2, A_{n-1})$}
in the taxonomy of \cite{Cecotti:2010fi}. 
The generalization mentioned in \S\ref{sec:other-ranks} similarly
corresponds to type $(A_{K-1}, A_{n-1})$.
The harmonic map equations are the equations giving Poincare
invariant vacua for this theory formulated on $S^1 \times \R^3$, 
with zero Wilson lines.

\subsection*{Acknowledgements} I am happy to 
thank Philip Boalch, David Dumas, Chris Fraser,
Laura Fredrickson and Michael Wolf
for very helpful discussions and explanations. 
I especially thank David Dumas for his assistance 
in understanding and using the code \cite{blaschke}, which was of
singular importance in confirming
the details of the picture described here. This work was supported
in part by NSF CAREER grant DMS-1151693 and a Simons Fellowship in Mathematics.

\section{Polynomials and polygons in \texorpdfstring{$\R\PP^2$}{RP2}} \label{sec:polynomials-and-polygons}

\subsection{A class of harmonic maps} \label{sec:harmonic-maps}

We consider harmonic maps 
\begin{equation}
 	g: \C \to \SO(3) \backslash \SL(3,\R).
\end{equation}
Given such a $g$, the corresponding \ti{Higgs field} is
\begin{equation}
	\varphi = - (\partial_z \tilde{g} \tilde{g}^{-1})^{+} \de z,
\end{equation}
where $A^+ = \half(A + A^t)$, and $\tilde{g}$ is any lift of $g$
to $\SL(3,\R)$.
The harmonicity of $g$ implies that $\varphi$ is holomorphic.
Thus if we define
\begin{equation}
  \phi_2 = -\frac12 \Tr \varphi^2, \qquad \phi_3 = - \frac13 \Tr \varphi^3,
\end{equation}
then $\phi_2$, $\phi_3$ are respectively holomorphic quadratic
and cubic differentials on $\C$.
In this paper we treat only the case of $g$ for which
\begin{equation} \label{eq:our-special-case}
  \phi_2 = 0, \qquad \phi_3 = P(z) \, \de z^3,
\end{equation}
for $P$ a polynomial of degree $n$, with complex coefficients.

The condition $\phi_2 = 0$ means that $g$ is not
only harmonic but also minimal \cite{MR2402597}.
For the purposes of this paper, this minimality is only indirectly
relevant: our main reason for focusing on the 
case \eqref{eq:our-special-case} 
is that it has been extensively studied in the recent \cite{MR3432157}.
We now briefly recall some facts established there; see
\cite{MR3432157} for a more detailed account and proofs.\footnote{In
comparing our formulas with those of \cite{MR3432157},
$\half C(z)_\there = P(z)_\here$.}

\begin{thm} For any polynomial $P$ there is a corresponding
harmonic $g$ for which \eqref{eq:our-special-case} holds.
This $g$ is unique up to the action of $\SL(3,\R)$ by isometries
on the target.
\end{thm}

The harmonic $g$ determines
an $(n+3)$-tuple of points $y_r \in \R\PP^2$, as follows.\footnote{This description
of the construction is different from that in \cite{MR3432157}, but produces
the same $y_r$.}
Write the leading term in $P(z)$ as $\mu z^n$ for some $\mu \in \C^\times$.
Fix a choice of root $\mu^{\frac{1}{n+3}}$. Then define $n+3$
rays $\ell_r \subset \C$ by
\begin{equation} \label{eq:antistokes-rays}
  \ell_r = \e^{\frac{2 \pi \I r}{n+3}} \mu^{-\frac{1}{n+3}} \R_+.
\end{equation}
For $z \in \ell_r$ we have $\mu z^{n+3} \in \R_+$,
so it admits a positive cube root, which we write
$(\mu z^{n+3})^\frac13$.
One of the eigenvalues of $\tilde g(z)$ is asymptotically
smallest,
behaving as
\begin{equation}
  \lambda \sim \exp \left(- \frac{3}{n+3} (\mu z^{n+3})^{\frac13} \right)
\end{equation} 
as $z \to \infty$ along $\ell_r$.
The corresponding eigenspace of $\tilde g(z)$
limits to $y_r \in \R\PP^2$.

\begin{thm} The $y_r$ are the vertices
of a convex $(n+3)$-gon in $\R\PP^2$. 
\end{thm}

Recall that $g$ is determined up to an overall $\SL(3,\R)$ action
from the right. This action transforms the $y_r$ by an overall $\SL(3,\R)$ action
on $\R\PP^2$. Thus, for a given polynomial $P$, the $y_r$ 
determine a convex $(n+3)$-gon 
in $\R\PP^2$ \ti{up to} this $\SL(3,\R)$ action.
Now:
\begin{itemize}
\item Let $\cM_{n+3} \subset (\R\PP^2)^{n+3} / \SL(3,\R)$
denote the moduli space of convex $(n+3)$-gons with vertices 
labeled. When $n > 0$, $\cM_{n+3}$ is a manifold,
diffeomorphic to $\R^{2n-2}$.

\item Let $\cP_n$ be the space of pairs $(P, \mu^{\frac{1}{n+3}})$
where $P$ is a degree $n$ polynomial and $\mu^{\frac{1}{n+3}}$ is
an $(n+3)$-rd root of the leading coefficient of $P$.
\end{itemize}
The passage from $(P, \mu^{\frac{1}{n+3}})$ to the $y_r$ defines a map
\begin{equation} \label{eq:map-Y}
 Y: \cP_n \to \cM_{n+3}.
\end{equation}
Both sides carry a natural action of $\Z / (n+3)\Z$;
the map $Y$ is equivariant for these actions.

Although we will not need this in what follows, we remark
that $Y$ is close to a homeomorphism, in the following sense.
$Y$ is invariant under affine-linear maps on $\cP_n$,
\begin{equation} \label{eq:affine-linear-symmetry}
 (P(z), \mu^{\frac{1}{n+3}}) \mapsto (P(a^{n+3} z+b), \mu^{\frac{1}{n+3}} a^n),  \qquad a \in \C^\times, b \in \C.
\end{equation}
Thus $Y$ descends to
\begin{equation}
 \overline{Y}: \overline\cP_n \to \cM_{n+3}
\end{equation}
where $\overline\cP_n$ is the quotient of $\cP_n$ by \eqref{eq:affine-linear-symmetry}.
For $n \ge 2$, it is shown in \cite{MR3432157} that $\overline Y$ is 
a homeomorphism.
A key building block of this result is a detailed theory of
error estimates for the harmonic map equation (Wang's equation), previously
developed in \cite{MR2144248,MR2350148,MR2402597}.

\subsection{The nonabelian Hodge correspondence} \label{sec:nonabelian-hodge}
In this section we explain how
$Y$ can be interpreted as a version of the \ti{nonabelian Hodge correspondence} between Higgs bundles and flat complex connections. 
This is not strictly necessary for reading the rest of the paper, 
but provides some context.

The nonabelian Hodge correspondence
originates from the celebrated work of Hitchin \cite{MR89a:32021}.
It was originally developed over a compact Riemann surface $C$ in
\cite{MR89a:32021,MR1179076,MR965220,MR887285},
and extended to the case of \ti{tame ramification} 
(first-order poles on $C$) in \cite{hbnc}.
We need the further extension to \ti{wild ramification}
(higher-order poles on $C$), developed in \cite{wnh}.

Briefly, 
the idea is as follows.
We consider the holomorphic bundle
$E = \cO^{\oplus 3}$ over $\C\PP^1$, equipped with the meromorphic Higgs field
\begin{equation}
  \varphi = \begin{pmatrix} 0 & 1 & 0 \\ 0 & 0 & 1 \\ -P(z) & 0 & 0 \end{pmatrix} \, \de z \enskip \in \enskip \End(E) \otimes K.
\end{equation}
The pair $(E,\varphi)$
is a Higgs bundle, with wild ramification at $z = \infty$.
More precisely, we consider
$(E,\varphi)$ as a \ti{good filtered Higgs bundle} in the language of \cite{MochizukiToda}, by assigning weights $(\alpha_1,\alpha_2,\alpha_3) = (\frac{n}{3}, 0, -\frac{n}{3})$
to the direct summands of $E$, and for a 
meromorphic section $s = (s_1,s_2,s_3)$ of $E$, 
defining $\nu_\infty(s) = \max \{ \nu_\infty(s_i) + \alpha_i \}_{i=1}^3$,
where $\nu_\infty(s_i) \in \Z$ is the ordinary order of singularity at $z = \infty$ of the
meromorphic function $s_i$.\footnote{The reason for the specific
weights $\alpha_i$ given here is that they arrange that
$\nu_\infty(\varphi(s)) / \de z = \frac{n}{3} + \nu_\infty(s)$
for any meromorphic section $s$ of $E$, which matches the 
singularity of the eigenvalues of $\varphi$,
$\lambda \sim c z^{\frac{n}{3}} \de z$.
More details of how
the nonabelian Hodge correspondence works in this example
will appear in \cite{FN}.}

We consider Hermitian metrics $h$ in $E$, 
compatible with the filtration: this means that for any meromorphic
section $s(z)$ of $E$, $h(s,s)$ scales like $\abs{z}^{2 \nu_\infty(s)}$ 
as $z \to \infty$. For any such $h$, let $D_h$ denote the Chern connection
in $E$, and $\varphi^{\dagger_h}$ the Hermitian adjoint of $\varphi$.
The key analytic fact which we use is existence and uniqueness for
\ti{Hitchin's equations}, established in this context
in \cite{wnh}:\footnote{More precisely, from the results of \cite{wnh} one
can directly deduce \autoref{thm:nonabelian-hodge} in the case $3 \vert n$.
The extension to other $n$ is believed to be straightforward, the basic mechanism being to pass to a cyclic covering, as described in \cite{MR1703088} --- see e.g. \cite{MochizukiToda} for the statement, though I do not know a reference
where all details have been explained.}
\begin{thm} \label{thm:nonabelian-hodge}
There exists a Hermitian metric $h$ in $E$, compatible with
the filtration, with
\begin{equation} \label{eq:hitchin-eq}
  F_{D_h} + [\varphi, \varphi^{\dagger_h}] = 0.
\end{equation}
This $h$ is unique up to overall scalar multiple.
\end{thm}
$h$ is the \ti{harmonic metric} associated to the Higgs bundle
$(E, \varphi)$.

We consider the family of complex connections in $E$ given by
\begin{equation} \label{eq:connection-family}
  \nabla(\zeta) = \zeta^{-1} \varphi + D_h + \zeta \varphi^{\dagger_h}, \qquad \zeta \in \C^\times.
\end{equation}
The equation \eqref{eq:hitchin-eq} implies that the 
connections $\nabla(\zeta)$ are all flat.
From the fact that $\varphi$ is traceless
it follows that the harmonic metric $h$ and flat connections
$\nabla(\zeta)$ are compatible with the standard volume form on $E$.
Thus each $\nabla(\zeta)$ is best thought of as a flat $\SL(3,\C)$-connection
over $\C\PP^1$, with an irregular singularity at $z = \infty$.
There is also an extra symmetry around: if we define a bilinear pairing
in $E$ by
\begin{equation}
  S = \begin{pmatrix} 0 & 0 & 1 \\ 0 & 1 & 0 \\ 1 & 0 & 0 \end{pmatrix} \enskip \in \enskip \Hom(E,E^*),
\end{equation}
then
\begin{equation}
  S^{-1} \varphi^T S = \varphi.
\end{equation}
This extra symmetry makes $(E,\varphi)$ into a \ti{real Higgs bundle} as
described in \cite{MR1174252}:
viewing $h$ as a map $E \to \overline{E}^*$, we have
the real structure $\tau = \overline{h}^{-1} \circ S: E \to \overline{E}$.
Then
\begin{equation}
  \overline{\nabla(\zeta)} = \overline{\zeta}^{-1} \overline{\varphi} + \overline{D_h} + \overline{\zeta} \overline{\varphi^{\dagger_h}} = \tau^{-1} \circ \nabla(\overline{\zeta}^{-1}) \circ \tau. 
\end{equation}
When $\abs{\zeta} = 1$ this becomes simply
\begin{equation}
  \overline{\nabla(\zeta)} = \tau^{-1} \circ \nabla(\zeta) \circ \tau,
\end{equation}
so using $\tau$ we can reduce $\nabla(\zeta)$ to an $\SL(3,\R)$-connection
in a real bundle $E_\R$.

The passage from the real Higgs bundle $(E,\varphi)$ to the flat $\SL(3,\R)$-connection
$\nabla(\zeta = 1)$ --- using the harmonic metric $h$ as intermediary --- 
is the nonabelian Hodge correspondence for real Higgs bundles.

Now we want to relate this to the harmonic maps $g$ described in
\S\ref{sec:harmonic-maps}.
This just involves a slight shift in perspective, following \cite{MR1179076,MR965220}.
Let $F$ denote the space of real flat sections
for $\nabla(\zeta = 1)$; $F$ is a $3$-dimensional real vector space
with a natural volume element.
Fix a basis $\{e_1, e_2, e_3\}$ of $F$, with unit volume.
Also choose a real $h$-unitary trivialization of the bundle $E$
away from $z = \infty$.
Then let
\begin{equation}
\tilde g(z) = (e_1(z), e_2(z), e_3(z)) \in \SL(3,\R).
 \end{equation}
Changing the basis of $F$ multiplies $\tilde g$ by a constant 
element of $\SL(3,\R)$ on the right;
changing the unitary trivialization of $E$ multiplies $\tilde g$ by
a smooth map $\C \to \SO(3)$ on the left.
Thus $\tilde{g}$ descends to a map $g: \C \to \SO(3) \backslash \SL(3,\R)$, determined
up to right-multiplication by an element of $\SL(3,\R)$. This is the
desired harmonic map.

Finally, we should explain how the polygons of \S\ref{sec:harmonic-maps}
arise. From our present point of view they have to do with the behavior
of the connection $\nabla(\zeta = 1)$ around the irregular singularity
at $z = \infty$. The rays $\ell_r$ defined in
\eqref{eq:antistokes-rays} are anti-Stokes rays.
Each $\ell_r$ determines a distinguished line in $F$,
consisting of flat sections with the fastest asymptotic decay
as $z \to \infty$: this gives the point $y_r$.

\section{Background for the integral iteration} \label{sec:background}

In this section we explain how to construct the input data needed in
the formulation of \autoref{conj:main}:
\begin{itemize}
\item parameters $R, \vartheta$ introduced in \S\ref{sec:rescaling} below,
\item a lattice $\Gamma$ equipped with a
\ti{period map} $Z$ and pairing $\IP{\cdot,\cdot}$,
described in \S\ref{sec:spectral-curve},
\item \ti{spectral coordinate} functions $X_\gamma$ on $\cM_{n+3}$,
described in \S\ref{sec:spectral-coordinates}-\S\ref{sec:asymptotic-abelianization-trees},
\item \ti{BPS counts} $\Omega(P_0, \gamma) \in \Z$, described in \S\ref{sec:bps-counts}.
\end{itemize}
In the process of constructing these data we will need the
notion of \ti{spectral network}, which we recall
in \S\ref{sec:spectral-network}-\S\ref{sec:network-asymptotics}.
The remaining \S\ref{sec:cluster} describes some relations between
our constructions and cluster algebra.

\subsection{Rescaling the cubic differential} \label{sec:rescaling}

From now on we take the polynomial $P$ of the form
\begin{equation} \label{eq:Pform}
  P(z) = R^3 \e^{-3 \I \vartheta} P_0(z), \qquad R \in \R_+, \ \vartheta \in \R.
\end{equation}
where $P_0(z)$ has only simple zeroes.

Of course, for a given $P(z)$ which we want to study, 
we could always set $R = 1$ and $\vartheta = 0$, by simply choosing
$P_0(z) = P(z)$. The point of our choice \eqref{eq:Pform} is:
\begin{itemize}
\item later we will want to study the asymptotic behavior as $R \to \infty$,
which is conveniently formulated in the parameterization \eqref{eq:Pform},
\item solving the problem for a single value of $\vartheta$ turns out to
give the solution for all $\vartheta$ at once, and so doing makes 
some of the structure more transparent.
\end{itemize}

\subsection{The spectral curve, homology and periods} \label{sec:spectral-curve}

Our first fundamental player is the \ti{spectral curve}, defined by
the equation
\begin{equation} \label{eq:spectral-curve}
	\Sigma = \{ x^3 + P_0(z) = 0 \} \subset \C^2.
\end{equation}
The projection
map $\pi: (x,z) \mapsto z$ makes
$\Sigma$ a branched $3$-fold cover of $\C$, with ramification
points of index $3$ over the $n$ zeroes of $P_0(z)$.
Using the Riemann-Hurwitz formula and looking at the ramification
around $z = \infty$ we compute
\begin{equation} \label{eq:sigma-topology}
  (\genus(\Sigma), \holes(\Sigma)) = \begin{cases} (n-2,3) & \text{ for } 3 \vert n, \\ 
  (n-1,1) & \text{ otherwise.} \end{cases}
\end{equation}
We will make frequent use of the lattice
\begin{equation}
  \Gamma = H_1(\Sigma, \Z),
\end{equation}
which has
\begin{equation}
  \rank(\Gamma) = 2n-2.
\end{equation}
Note that this formula is uniform in $n$, despite the case structure
in \eqref{eq:sigma-topology}. 
Also note a numerical ``coincidence'' which will be important later:
$\rank(\Gamma) = \dim \cM_{n+3}$.

$\Gamma$ is equipped with the skew-symmetric
intersection pairing $\IP{\cdot,\cdot}$, and the
\ti{period homomorphism}
\begin{equation} \label{eq:periods}
	Z: \Gamma \to \C, \qquad   Z_\gamma = \oint_\gamma x \, \de z.
\end{equation}

\subsection{The WKB spectral network} \label{sec:spectral-network}

Given $(P_0, \vartheta)$ there is a corresponding \ti{WKB spectral network}
$\cW(P_0, \vartheta)$. 
Two examples of WKB spectral networks $\cW(P_0, \vartheta)$ --- the
only two which we will consider in detail in this paper --- 
are shown 
in Figures \ref{fig:pentagon-network} and \ref{fig:hexagon-network} below.

In the rest of this section we describe what $\cW(P_0, \vartheta)$ is
and how it is constructed. Our description uses the language of 
\cite{Gaiotto2012};
essentially the same networks had been discovered earlier as
\ti{Stokes graphs}, e.g.
\cite{berk:988,MR2132714}.\footnote{The networks described here may look
surprising to a reader familiar with \cite{Gaiotto2012,berk:988,MR2132714}; in those references the initial trajectories emanate 
from $3$-valent vertices, while here we have $8$-valent vertices. 
The reason
for this difference is that we are studying a situation where all
ramification points of $\Sigma$ have index $3$, while those earlier 
references mostly concerned ramification points of index $2$. Had we chosen
to study general harmonic maps as opposed to minimal ones, we would
generically get ramification points of index $2$, and our pictures
would look more like those of \cite{Gaiotto2012,berk:988,MR2132714}.
The relation between these two situations is illustrated in
Figure 36 of \cite{Gaiotto2012}. For a quick definition of the WKB
spectral network in the more generic situation see Section 4.2
of \cite{cluster-susy}.}

The WKB spectral network is a collection of \ti{WKB $\vartheta$-trajectories}.
A WKB $\vartheta$-trajectory is a path $z(t)$ on $\C$, obeying 
a first-order ODE, depending on a choice
of an ordered pair of distinct sheets
of $\Sigma$ lying over $z(t)$. We label these sheets as
$(x_i(t), x_j(t))$ or sometimes just $(i,j)$.
The equation is:
\begin{equation} \label{eq:trajectory-ode}
  (x_i(t) - x_j(t)) \frac{\de z(t)}{\de t} = \e^{\I \vartheta}.
\end{equation}
Note that \eqref{eq:trajectory-ode} is invariant
under the operation of simultaneously reversing $i \leftrightarrow j$ and $t \leftrightarrow -t$.

On any simply connected patch $U \subset \C$ with $P_0(z) \neq 0$ on $U$,
the WKB $\vartheta$-trajectories make up three foliations, labeled by the three
choices of unordered pair $(i,j)$; the three foliation leaves passing through
each point meet at angles $2 \pi/3$.
Around a simple zero $z_0$ of $P_0(z)$ the structure is more interesting:
there are $8$ backward inextendible WKB $\vartheta$-trajectories which end on $z_0$; call these \ti{critical trajectories}. See Figure \ref{fig:near-branch}.
\insfig{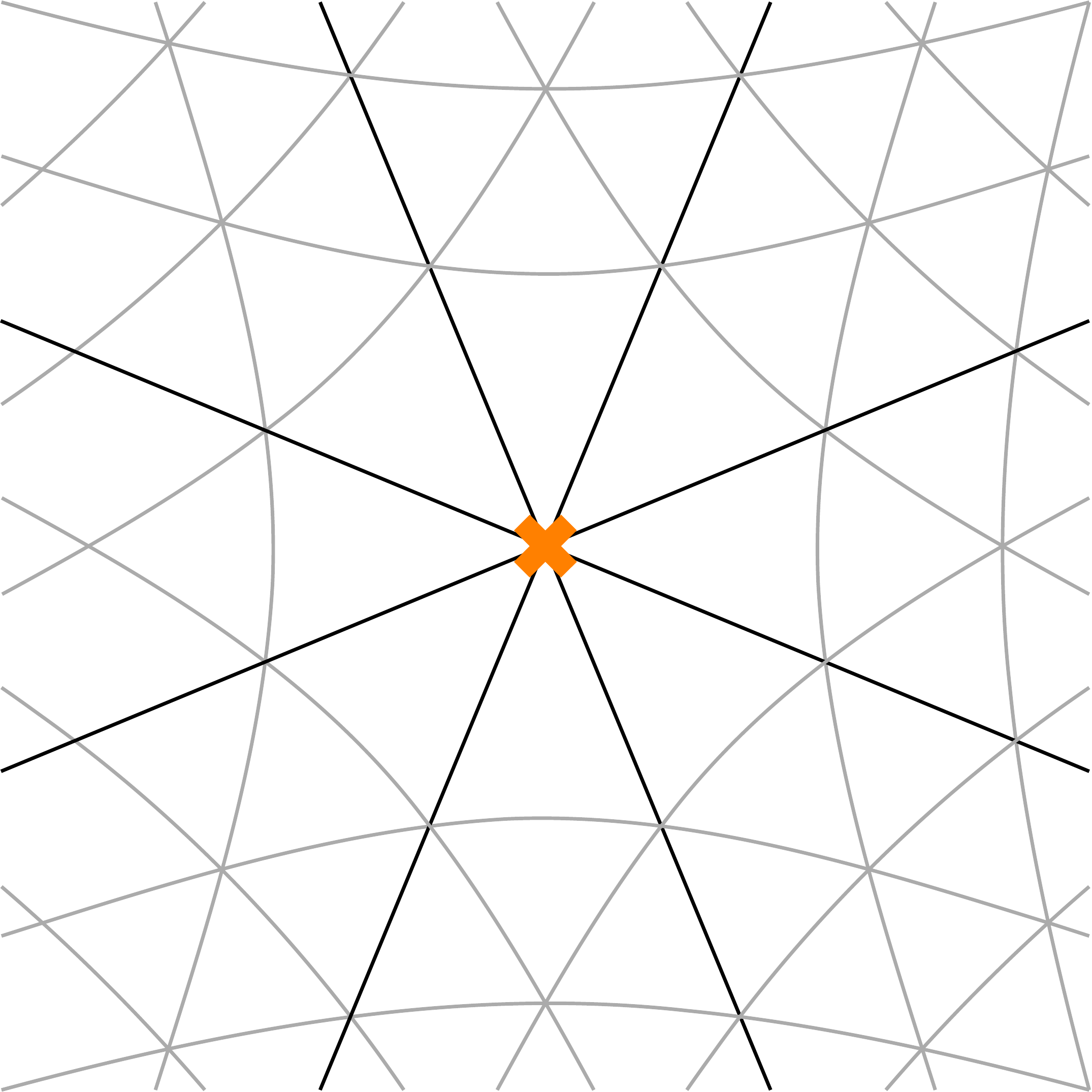}{0.18}{Some WKB $\vartheta$-trajectories in a neighborhood of a simple zero of $P_0$ (orange cross).
The critical trajectories are shown in black, others in gray.}

The network $\cW(P_0, \vartheta)$ is constructed as follows.
We begin with the $8n$ critical trajectories emanating from the
$n$ zeroes of $P_0$, and extend them to $t \to +\infty$ by integrating
the ODE \eqref{eq:trajectory-ode}. These trajectories will be included
in $\cW(P_0, \vartheta)$.
Next we iteratively add more trajectories to $\cW(P_0, \vartheta)$, 
as follows.

We consider intersections between
trajectories already included in $\cW(P_0, \vartheta)$.
For each intersection there are three possibilities: either 
the trajectories meet head-on
(in which case they actually coincide, differing only in their
orientation and reversal of the sheet labels),
they intersect in an angle $\frac{\pi}{3}$, or they intersect in 
an angle $\frac{2 \pi}{3}$.
For each intersection in angle $\frac{2\pi}{3}$ we add a new
trajectory, as follows.
The fact that the intersection angle is $\frac{2 \pi}{3}$
implies that the labels of the intersecting trajectories are of the form
$(i_1,j_1) = (i,j)$ and $(i_2,j_2) = (j,k)$.
We add a new WKB $\vartheta$-trajectory beginning from the
intersection point, 
with the label $(i, k)$, as shown in Figure
\ref{fig:ijkl-collision-birth}.

\insfig{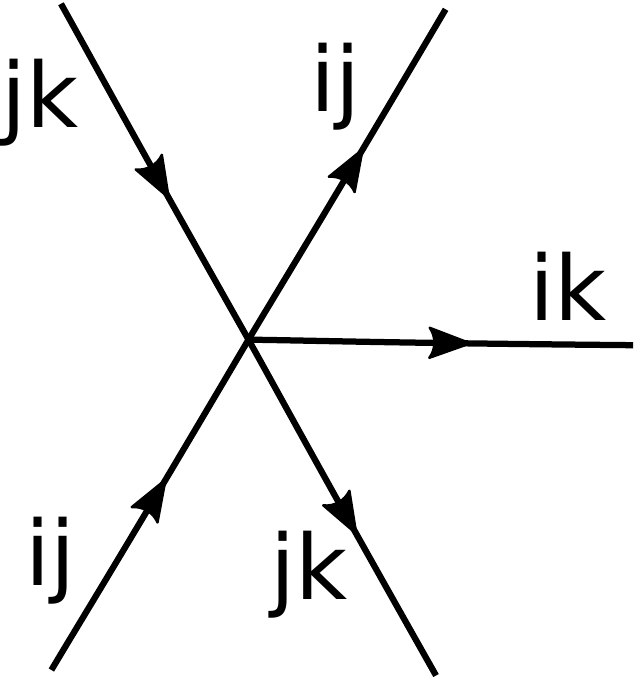}{0.33}{When two $\vartheta$-trajectories
in $\cW(\vartheta, P_0)$ with labels $(i_1,j_1) = (i,j)$ 
and $(i_2,j_2) = (j,k)$ intersect, we add to $\cW(\vartheta, P_0)$ 
a new $\vartheta$-trajectory carrying labels $(i_3,j_3) = (i,k)$,
beginning at the intersection point.}

As before, we extend these new trajectories to $t \to +\infty$. This may create
new intersections between trajectories. We then repeat the process,
letting these new intersections give birth to new trajectories,
and so on. Define $\cW(P_0, \vartheta)$ to be the full collection
of trajectories produced in this fashion. 
\ti{A priori} there is no reason this collection should be 
finite, but in the examples we study explicitly 
in \S\ref{sec:case-studies} below, it is finite.

There is one phenomenon which requires extra attention: when we try
to extend a trajectory in $\cW(P_0, \vartheta)$ to 
$t \to + \infty$, it may run into a zero of $P_0$ at some finite $t$.
In this case we say that $(P_0, \vartheta)$ is
\ti{BPS-ful}.
The BPS-ful case is also the case in which there is at least one 
head-on collision between trajectories in $\cW(P_0, \vartheta)$.
If this does not happen then
we say $(P_0, \vartheta)$ is \ti{BPS-free}.

\subsection{Asymptotics of WKB spectral networks} \label{sec:network-asymptotics}

The behavior of $\cW(P_0, \vartheta)$ near $z \to \infty$ is particularly
simple: the WKB $\vartheta$-trajectories approach $2n+6$ asymptotic 
directions. It is convenient to compactify to
 $\overline\C = \C \sqcup S^1$, introducing 
a ``circle at infinity'' (i.e. $\overline\C$ is the oriented real blow-up
of $\C\PP^1$ at $z = \infty$). The $2n+6$ asymptotic
directions then give $2n+6$ marked points on this circle.

Each marked point is labeled by an ordered pair of sheets $(i,j)$, 
giving the label for all the WKB $\vartheta$-trajectories asymptotic 
to it.
When we move from one marked point to the next
(going around say counterclockwise)
one label stays the same while the other changes, alternating between
first and last: i.e. the labels on consecutive rays follow 
the pattern $ij$, $ik$, $jk$, $ji$, $ki$, $kj$, \dots. 
Again see Figure \ref{fig:pentagon-network} and
Figure \ref{fig:hexagon-network} for illustrative examples.\footnote{
In comparing the asymptotics in those figures to the description 
above, one must keep in mind the permutations
of sheets which occur when one crosses the branch cuts. For example,
in Figure \ref{fig:pentagon-network}, starting from the rightmost marked
point just above the branch cut, we see the sequence
$31,21,23,13,12$; the next label would ordinarily be $32$, 
but we also cross
a branch cut which induces the permutation $(123)$, so the next
label is instead $13$.}

These marked points divide the circle at infinity into $2n+6$ arcs.
Call an arc at infinity
\ti{initial} (\ti{final}) if its two boundary points have the same 
initial (final) label; there are $n+3$ arcs of each type. 
Each asymptotic direction $\ell_r$ lies at the midpoint of
one of the final arcs.
Thus each $\ell_r$ can be labeled by the final label for the
two nearest boundary rays; call this the \ti{fading sheet}
at $\ell_r$.

This asymptotic structure is essentially universal: changing
$(P_0, \vartheta)$ changes it only by an overall rotation
of the circle at infinity.

\subsection{Spectral coordinates} \label{sec:spectral-coordinates}

One of the predictions of \cite{Gaiotto:2008cd,Gaiotto:2009hg,Gaiotto2012}
is that any BPS-free $(P_0, \vartheta)$ determines a coordinate system on 
$\cM_{n+3}$.\footnote{Strictly speaking, what \cite{Gaiotto:2008cd,Gaiotto:2009hg,Gaiotto2012} predict is a \ti{local}
coordinate system on a patch of a \ti{complexification} of $\cM_{n+3}$,
but it seems
reasonable to conjecture that we get
\ti{global} coordinates after restricting to $\cM_{n+3}$.}
More specifically, $(P_0,\vartheta)$ 
should determine for each $\gamma \in \Gamma$
a function $X_\gamma: \cM_{n+3} \to \R^\times$, such that
\begin{equation}
	X_{\gamma + \mu} = X_\gamma X_\mu,
\end{equation}
and if we choose generators $\{\gamma_i\}_{i=1}^{2n-2}$ for $\Gamma$,
the functions $\{X_{\gamma_i}\}_{i=1}^{2n-2}$ should give a coordinate system
on $\cM_{n+3}$.
When $(P_0, \vartheta)$ is varied while remaining BPS-free,
the $X_\gamma$ should not change; when $(P_0, \vartheta)$ is varied across
a BPS-ful locus, the $X_\gamma$ may change.

\subsection{Abelianization trees} \label{sec:abelianization-trees}

In this section and the next, 
we describe how the spectral coordinates $X_\gamma$
are constructed using the spectral networks
$\cW(P_0, \vartheta)$.

Define an \ti{abelianization tree compatible with $\cW(P_0,\vartheta)$}
to be a collection of oriented arcs in $\overline\C$, with
each arc labeled by a sheet of $\Sigma$ and a representation
of $\SL(3,\R)$ (either fundamental $V$ or its dual $V^*$),
with the following properties:
\begin{itemize}
  \item Each arc has two endpoints. The initial point of each
  arc may lie at one of the $\ell_r$ on the circle at infinity, or 
  else lie at a junction as shown in 
  Figure \ref{fig:junctions}. The endpoint of each arc
  must lie at a junction.
  \insfig{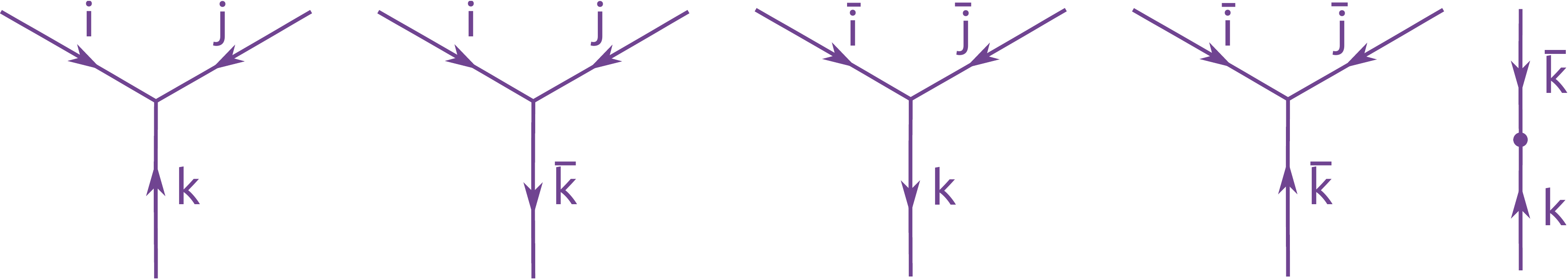}{0.36}{The types of junctions allowed in an abelianization
  tree. For the trivalent junctions $i,j,k$ must all be distinct. If the
  sheet label carries an overbar, then the arc is carrying 
  representation $V^*$, otherwise it is carrying representation $V$.}
  \item If an arc ends at $\ell_r$, then that arc carries the
  representation $V$, and its sheet label matches the fading sheet
  at $\ell_r$ (as defined in \S\ref{sec:network-asymptotics}).
  \item Arcs of the abelianization tree do not cross one another.
  \item No arc carrying the label $i$ and representation $V$ 
  crosses a trajectory of $\cW(P_0,\vartheta)$ carrying a label $(i,j)$.
  \item No arc carrying the label $i$ and representation $V^*$ 
  crosses a trajectory of $\cW(P_0,\vartheta)$ carrying a label $(j,i)$.
\end{itemize}
Some examples of abelianization trees compatible with spectral
networks appear in Figure 
\ref{fig:pentagon-gamma1-trees} and Figure 
\ref{fig:hexagon-gamma1-trees} below.

Dropping the sheet labels from an abelianization tree $h$ induces
a \ti{tensor diagram} on $\overline\C$.
This is a notion with a
long history: see \cite{Fomin2012} for a very clear and precise review.
This tensor diagram determines an $\SL(3,\R)$-invariant map
\begin{equation}
 A_h: V^{n+3} \to \R,
\end{equation}
using the standard intertwiners 
\begin{equation}
  V \otimes V \otimes V \to \R, \quad
V \otimes V \to V^*, \quad 
V^* \otimes V^* \to V, \quad
V^* \otimes V^* \otimes V^* \to \R, \quad 
V \otimes V^* \to \R.
\end{equation}
(We use the standard orientation
of $\C$ to fix the orderings where needed.)
$A_h$ is homogeneous: its scaling weights under $(\R^\times)^{n+3}$
are the numbers of arcs ending at the $n+3$ marked points $\ell_r$.

\subsection{Combining abelianization trees for spectral coordinates} \label{sec:combining-abelianization-trees}

Now fix a formal linear combination of abelianization trees
$\sum w_m h_m$, with weights $w_m \in \Z$, such that the total
weighted number of arcs ending at each point $\ell_r$ is zero.
Then the function
\begin{equation}
  X = \prod_m A_{h_m}^{w_m}
\end{equation}
is invariant under $(\R^\times)^{n+3}$, and so descends
to $\cM_{n+3} \subset (\R\PP^2)^{n+3} / \SL(3,\R)$.
$X$ will be one of the
spectral coordinates $X_\gamma$, for some
$\gamma \in \Gamma$. It only remains to explain what
$\gamma$ is.

For this purpose, note each arc $p$ of the abelianization tree $h_m$ 
has a canonical lift to a $1$-chain $p^\Sigma$
on $\Sigma$: namely, if $p$ is labeled by the sheet $i$, $p^\Sigma$
is the lift of $p$ to sheet $i$, with the orientation
of $p$ if $p$ is carrying representation $V$, otherwise with
the opposite orientation.
Summing these lifts gives a singular $1$-chain $h_m^\Sigma$ for each 
abelianization tree $h_m$.

We now consider the $1$-chain
\begin{equation}
  c = \sum_m w_m h_m^\Sigma.
\end{equation}
$c$ is generally not closed,
because of the trivalent
junctions. However, the boundary $\partial c$ is a $0$-chain 
pulled back from the base $\overline \C$, so $c$
becomes closed if we project to the quotient
of the singular chain complex $C_*(\Sigma)$ by the subgroup of chains
pulled back from $\overline \C$. Moreover 
this projection is an isomorphism on homology, because 
$\overline \C$ is contractible.
Thus $c$ determines a class
\begin{equation}
  \left[ \sum_m w_m h_m^\Sigma \right] \in \Gamma.
\end{equation}
Now we can finally state our definition of the spectral coordinate:
\begin{equation} \label{eq:def-spectral}
  X_{\left[ \sum_m w_m h_m^\Sigma \right]} = \prod_m A_{h_m}^{w_m}.
\end{equation}

\eqref{eq:def-spectral} defines $X_\gamma$ for any $\gamma \in \Gamma$
which can be realized as $\gamma = \left[ \sum_m w_m h_m^\Sigma \right]$
for some abelianization trees $h_m$ and weights $w_m$.
For this definition to be unambiguous, it must be true that
every $\gamma \in \Gamma$ admits exactly one
such realization. 
This is true in the examples we consider in \S\ref{sec:case-studies} 
below. I do not know how generally it holds; in 
cases where it fails, we would have to resort to a more complicated
construction of the $X_\gamma$, beginning from the notion of
\ti{abelianization} described in \cite{Gaiotto2012} and further reviewed 
in \cite{cluster-susy}.

\subsection{Asymptotic abelianization trees} \label{sec:asymptotic-abelianization-trees}

For different choices of $(P_0, \vartheta)$ 
we get different compatible abelianization trees, and thus 
different spectral coordinates $X_\gamma$. 
However, as we have noted in \S\ref{sec:network-asymptotics}, 
the structure of $\cW(P_0, \vartheta)$ near $z = \infty$ is
independent of $(P_0, \vartheta)$, up to an overall rotation of the
plane. Thus any abelianization tree which lies entirely in the
asymptotic region is in a sense universal, existing for
every $(P_0, \vartheta)$. We call these \ti{asymptotic abelianization trees}.

\insfig{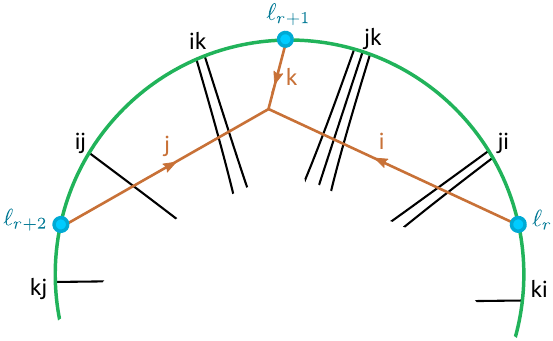}{1.5}{An asymptotic abelianization tree.
Black lines represent WKB $\vartheta$-trajectories in the asymptotic region
of a WKB spectral network.}

In Figure \ref{fig:asymptotic-abelianization-tree} we show an 
example of an asymptotic abelianization tree $h$, depending
on a choice of $r \in \{ 1, \dots, n+3 \}$. 
The corresponding
holonomy function is
\begin{equation} \label{eq:consecutive-plucker}
  A_h = p(r,r+1,r+2)
\end{equation}
where $p(a,b,c)$ is a
a ``Pl\"ucker coordinate'' depending on
$3$ of the $v_i$,
\begin{equation} \label{eq:plucker}
  p(a,b,c): (v_i)_{i=1}^{n+3} \mapsto \det( v_a, v_b, v_c).
\end{equation}

\subsection{BPS counts} \label{sec:bps-counts}

The last ingredient we need is a collection of
\ti{BPS counts} $\Omega(P_0, \gamma) \in \Z$ for $\gamma \in \Gamma$.
Here we describe one geometric approach to defining these counts.
(This method is practicable for the simple examples
we treat in \S\ref{sec:case-studies}, but for more elaborate
examples it would be hard to use in practice. Thus
we should mention that there are other ways,
e.g. the \ti{mutation method} 
\cite{Gaiotto:2010be,Alim:2011kw,Alim:2011ae} or the
\ti{spectrum generator / DT transformation} 
\cite{Gaiotto:2009hg,Goncharov2016,Longhi2016}.)

Any WKB $\vartheta$-trajectory $p$ admits a canonical lift
to a $1$-chain $p^\Sigma$ on $\Sigma$: namely, if $p$ is labeled by the
pair of sheets $(i,j)$, $p^\Sigma$ is the closure of the lift of $p$ to 
sheet $i$, plus the orientation-reversal of the 
closure of the lift of $p$ to sheet $j$.
For $\gamma \in \Gamma$ a \ti{finite web of charge $\gamma$} 
is a collection of trajectories $p_l$ such that $\sum_l p_l^\Sigma$
is a compact $1$-cycle in homology class $\gamma$.
From the equation \eqref{eq:trajectory-ode}
it follows that if there is a finite web of charge $\gamma$
inside a WKB spectral network $\cW(P_0, \vartheta)$, then we must have
$\vartheta = \arg Z_\gamma$. Moreover, in this case the 
network $\cW(P_0, \vartheta)$ is BPS-ful in the sense of \S\ref{sec:spectral-network}.

The desired $\Omega(P_0, \gamma)$ is a count
of the finite webs of charge $\gamma$ which occur
inside the network $\cW(P_0, \vartheta = \arg Z_\gamma)$.
The meaning of the word ``count'' here is given by an algorithm
explained in \cite{Gaiotto2012}. We will not describe this algorithm
here, since in full generality it is a bit complicated; 
for the examples of \S\ref{sec:case-studies} we 
only need some simple
special cases, where the finite web is either a single critical trajectory
connecting two zeroes of $P_0$ or a three-string junction: 
see Figure \ref{fig:finite-webs}.
In either of these cases, the algorithm in \cite{Gaiotto2012}
says that the finite web contributes $+1$ to $\Omega(P_0, \gamma)$.
\insfig{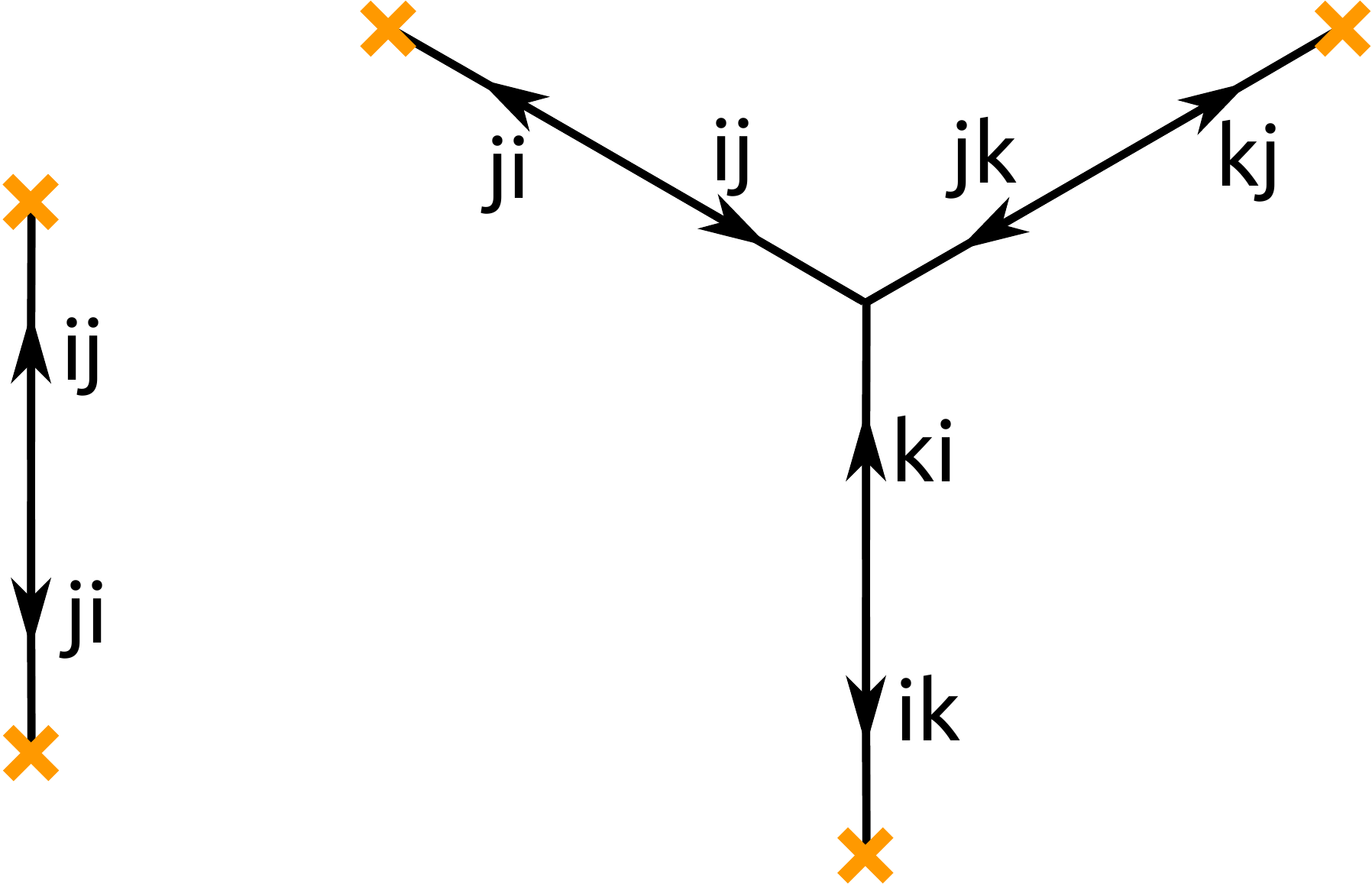}{0.3}{Two examples of finite webs.}

A few final remarks about the BPS counts:

\begin{itemize}
\item For the examples of $P_0$ we consider in \S\ref{sec:case-studies}, there
exist only finitely many $\gamma \in \Gamma$ for which
$\Omega(P_0, \gamma) \neq 0$. We conjecture that this is always
true for $n \le 5$.

\item
The BPS counts have the symmetry property
\begin{equation} \label{eq:Omega-symmetry}
  \Omega(P_0, \gamma) = \Omega(P_0, -\gamma).
\end{equation}
To see this, note that $\cW(P_0, \vartheta)$ and $\cW(P_0, \vartheta+\pi)$
differ only by relabeling all trajectories $ij \to ji$, and thus
any finite web in $\cW(P_0, \vartheta)$ with charge $\gamma$
has a partner in $\cW(P_0, \vartheta+\pi)$ with charge $-\gamma$.

\item
Similarly, the $\Omega(P_0, \gamma)$ are invariant under the
$\Z/3\Z$ action on $\Gamma$ induced by the $\Z/3\Z$ action on $\Sigma$
(cyclic permutation of sheets); this reflects the fact that
$\cW(P_0, \vartheta)$ and $\cW(P_0, \vartheta + \frac{2\pi}{3})$
differ only by cyclic permutation of the sheet labels.

\item 
The $\Omega(P_0, \gamma)$ are expected
to be examples of \ti{generalized Donaldson-Thomas invariants} in the sense of
Kontsevich-Soibelman / Joyce-Song \cite{ks1,Joyce:2008pc}.
In particular, the collection $\{\Omega(P_0, \gamma)\}_{\gamma \in \Gamma}$
depends on $P_0$ in a piecewise constant fashion;
it can jump where $P_0$ crosses a \ti{wall of marginal stability} in the
space of polynomials. These jumps are governed by the Kontsevich-Soibelman
wall-crossing formula \cite{ks1,Gaiotto:2008cd,Gaiotto2012}.

\item
A related fact is that, when $(P_0, \vartheta)$ cross a BPS-ful locus,
the jump of the spectral coordinates $X_\gamma$ is determined by
the $\Omega(P_0, \gamma)$ for those $\gamma$ such that
$Z_\gamma \e^{-\I \vartheta} \in \R_-$.
\end{itemize}

\subsection{Spectral coordinates as cluster coordinates} \label{sec:cluster}

This section describes an interpretation of the functions
$A_h$ associated to abelianization trees in the language of cluster algebra.
It is not strictly necessary for the rest of the paper.
See \cite{Fomin2012} for much more on the cluster structures
we use here, and their relation to tensor diagrams.

Let $\cC_{n+3}$ denote the coordinate ring of $V^{n+3} / \SL(3,\R)$.
This ring is a classic example of a \ti{cluster algebra},
as first described in \cite{MR2205721}.
This structure picks out a collection of distinguished \ti{clusters}
in $\cC_{n+3}$: these are collections of \ti{cluster $\fA$-variables} $\{a_m\}_{m=1}^{3n+1}$,
where the first $2n-2$ are called \ti{unfrozen} and depend on 
the choice of cluster,
and the last $n+3$ are \ti{frozen} and are common to all clusters.

The reason for discussing the clusters here is the 
following: there seems to be a correspondence between
clusters and spectral networks (as anticipated in \cite{Gaiotto2012}.) Namely, given a spectral network $\cW$,
the functions $A_h$ realized by abelianization trees compatible with $\cW$
are often the monomials in the cluster $\fA$-variables
of some cluster in $\cC_{n+3}$. This is true for every spectral
network $\cW(P_0,\vartheta)$ I have examined (with $n = 2,3$).
I conjecture that it is true for all spectral
networks $\cW(P_0,\vartheta)$ for BPS-free $(P_0, \vartheta)$
when $n \le 5$.\footnote{These 
are the cases for which $\cC_{n+3}$ has only 
finitely many clusters \cite{MR2205721}.}
It would be very interesting to understand the relation
between this proposal and the results in \cite{Fomin2012},
and to develop an efficient method of mapping the spectral networks
to the clusters.

In this correspondence, the $n+3$ asymptotic abelianization trees described 
in \S\ref{sec:asymptotic-abelianization-trees} give the $n+3$ frozen variables,
which are the Pl\"ucker coordinates \eqref{eq:consecutive-plucker}.
In the language of \cite{Fomin2012} these are \ti{special
invariants}.
Non-asymptotic abelianization trees give the $2n-2$ unfrozen variables.
Some of these unfrozen variables are also Pl\"ucker coordinates $p(a,b,c)$,
now with $a,b,c$ not consecutive.
For example, when $n=2$, all cluster $\fA$-variables
are Pl\"ucker coordinates.
The first example of a cluster $\fA$-variable
which is not a Pl\"ucker coordinate arises for $n=3$; we will encounter
it in \S\ref{sec:hexagon} below (see \eqref{eq:hexapod-invariant}).

The cluster structures we discussed above are closely related to ones
introduced in \cite{MR2233852}, which were
very important in the original developments
leading to \cite{Gaiotto2012} and \autoref{conj:main}.
However, as far as I understand, the cluster structures we are now
considering are not \ti{quite} examples of the formalism 
in \cite{MR2233852}: that formalism does include
moduli spaces of flat connections
with irregular singularity, but a different sort
of irregular singularity than we need here,
which leads to moduli spaces involving complete
flags rather than lines.

\section{The integral iteration and its consequences} \label{sec:integral-and-consequences}

We can now give the sharp formulation of
\autoref{conj:main} and some consequences.

\subsection{The integral iteration} \label{sec:integral-equations}

In this section we formulate \autoref{conj:main} precisely. 
Concretely, we give a conjectural way of computing the
point $Y(P) \in \cM_{n+3}$,
in the spectral coordinate system $X_\gamma$ on $\cM_{n+3}$
associated to $(P_0, \vartheta)$.\footnote{Of course, once we have 
determined the coordinates of $Y(P)$ in one coordinate
system on $\cM_{n+3}$ we can then change coordinates to any other; nevertheless, \autoref{conj:main} is most naturally phrased
as a recipe which produces specifically the spectral
coordinates $X_\gamma$.}

Although our aim is to compute \ti{numbers} $X_\gamma$, the strategy
is first to construct \ti{functions} $\cX_\gamma(\zeta)$ of a parameter
$\zeta \in \C^\times$, as follows.
We begin with $\cX_\gamma^{(0)}(\zeta) = 0$
and then define inductively\footnote{In \eqref{eq:iteration} the factor
$\log(1 + \cX_\mu(\zeta'))$ appears. In the general formalism of
\cite{Gaiotto:2008cd} one would expect to see instead
$\log(1 - \sigma(\mu) \cX_\mu(\zeta'))$ for some
$\sigma: \Gamma \to \{ \pm 1 \}$. In the examples of this paper we always
have $\sigma(\mu) = -1$ whenever $\Omega(P_0,\mu) \neq 0$, so we
have simplified by making this substitution.}
\begin{multline} \label{eq:iteration}
\cX_\gamma^{(k+1)}(\zeta) = \exp \Bigg[ R(\zeta^{-1} Z_\gamma + \zeta \overline{Z}_\gamma) \\ + \frac{1}{4
\pi \I} \sum_{\mu \in \Gamma} \Omega(P_0, \mu) \langle \gamma,\mu
\rangle \int_{Z_\mu \R_-} \frac{\de \zeta'}{\zeta'} \frac{\zeta' +
\zeta}{\zeta' - \zeta} \log (1 + \cX^{(k)}_{\mu}(\zeta'))\Bigg].
\end{multline}
The desired functions are obtained in the limit:
\begin{equation} \label{eq:X-limit}
  \cX_\gamma(\zeta) = \lim_{k \to \infty} \cX^{(k)}_\gamma(\zeta).
\end{equation}
A sketch of an argument for the existence of this limit at sufficiently
large $R$
can be found in \cite{Gaiotto:2008cd}. For our purposes in this paper
we take the convergence as a working assumption.
Finally, to recover the desired
numbers $X_\gamma$, we specialize:
\begin{equation} \label{eq:X-specialize}
  X_\gamma = \cX_\gamma(\zeta = \e^{\I \vartheta}).
\end{equation}
The equations \eqref{eq:iteration}, \eqref{eq:X-limit}, 
\eqref{eq:X-specialize}, together with the
definitions of symbols given in \S\ref{sec:background},
make up the precise form of \autoref{conj:main}.

A few remarks:
\begin{itemize}
\item The right side of \eqref{eq:X-specialize} 
is real, as it should be,
while $\cX_\gamma(\zeta)$ is generally complex for $\abs{\zeta} \neq 1$.
The proof that $\cX_\gamma(\zeta = \e^{\I \vartheta})$ 
is real uses the symmetry property \eqref{eq:Omega-symmetry}.
\item The $\cX_\gamma(\zeta)$ depend on $\zeta$ in a \ti{piecewise}
holomorphic fashion, jumping whenever $\zeta$ crosses the 
contours of integration in \eqref{eq:iteration}, i.e. the
rays $Z_\mu \R_-$ with $\Omega(P_0, \mu) \neq 0$ and
$\IP{\mu,\gamma} \neq 0$. This means that
the $X_\gamma$ we compute from \eqref{eq:X-specialize}
depend discontinuously on $\vartheta$. Of course, the polygon
$Y(P) \in \cM_{n+3}$ depends continuously on $\vartheta$. The 
jumps of $X_\gamma$ just reflect
the fact that the spectral coordinate system on $\cM_{n+3}$
which we use jumps when $\vartheta$ crosses a BPS-ful phase.
\item The $\cX_\gamma(\zeta)$ for general $\zeta$ 
do have an interpretation,
most easily expressed in the language of \S\ref{sec:nonabelian-hodge}:
they are the spectral coordinates of the family
of connections \eqref{eq:connection-family}.
\end{itemize}

\subsection{Exact consequences} \label{sec:exact}

Now let us describe some consequences of \autoref{conj:main}.

One important special case arises
when $\gamma$ lies in the kernel of the pairing $\IP{\cdot,\cdot}$.
(Recall from \S\ref{sec:spectral-curve} that this 
kernel is nontrivial only when $3 \vert n$.)
In this case the coefficients of the integrals in \eqref{eq:iteration}
vanish, and \eqref{eq:X-limit} becomes the exact prediction
\begin{equation}
  \cX_\gamma(\zeta) = \exp \left[ R(\zeta^{-1} Z_\gamma + \zeta \overline{Z}_\gamma) \right],
\end{equation}
giving after the substitution \eqref{eq:X-specialize}
\begin{equation} \label{eq:exact-specialized}
  X_\gamma = \exp\left(a_\gamma R\right),
\end{equation}
with
\begin{equation} \label{eq:a-gamma}
  a_\gamma = 2 \re(\e^{-\I \vartheta} Z_\gamma).
\end{equation}

\subsection{Asymptotic consequences for \texorpdfstring{$R \to \infty$}{R -> infinity}} \label{sec:asymptotics}

The more interesting case is that of $\gamma$
\ti{not} in the kernel of the pairing $\IP{\cdot,\cdot}$.
Then \eqref{eq:X-limit} does not reduce to a simple 
exact formula for the spectral coordinate $X_\gamma$.
Nevertheless, we can extract asymptotic formulas in
the limit $R \to \infty$. For this purpose define
\begin{equation}
  \alpha_\mu = -\frac{Z_\mu}{\abs{Z_\mu}}  
\end{equation}
and assume for simplicity that, for all $\mu \in \Gamma$
such that $\Omega(P_0, \mu) \neq 0$, we have $\zeta \neq \alpha_\mu$.
We consider the integral equation obeyed 
by $\cX_\gamma(\zeta)$ (obtained by setting $k = \infty$ in 
\eqref{eq:iteration}) and make a self-consistent analysis:
assume that the integral terms are exponentially small in the 
$R \to \infty$ limit, use this assumption to replace $\cX_\mu(\zeta')$ 
by $\exp R(\zeta'^{-1} Z_\mu + \zeta' \overline{Z}_\mu)$ in the integral,
then evaluate the $R \to \infty$ asymptotics of the integral 
by the saddle point method.

This leads to the following prediction.
Write
\begin{equation} \label{eq:first-asymptotic-correction}
\cX_\gamma(\zeta) = \exp \left[ (\zeta^{-1} Z_\gamma + \zeta \overline{Z}_\gamma)R + \sum_{\mu \in \Gamma} \frac{\Omega(P_0, \mu) \IP{\gamma,\mu}}{4 \pi \I}
  \frac{\alpha_\mu + \zeta}{\alpha_\mu - \zeta} \sqrt{\frac{\pi}{R \abs{Z_\mu}}} \e^{-2 \abs{Z_\mu} R} + \delta(R,\zeta) \right],
\end{equation}
for some ``remainder'' $\delta(R,\zeta)$.\footnote{The values $\delta(R,\zeta) = \pm \infty$ are allowed, but our claim \eqref{eq:asymptotic-claim} implies that there is some $R_0$ such that $\delta(R,\zeta)$ is finite for $R > R_0$.}
We predict that
\begin{equation} \label{eq:asymptotic-claim}
\lim_{R \to \infty} \sqrt{R} \e^{2 \rho R} \delta(R,\zeta) = 0,
\end{equation}
where
\begin{equation}
  \rho = \min \left\{ \abs{Z_\mu}: \mu \in \Gamma, \ \Omega(P_0, \mu) \IP{\gamma,\mu} \neq 0 \right\}.
\end{equation}
Said otherwise: as $R \to \infty$, the leading behavior of $\log \cX_\gamma$ is of order $R$;
the first correction is exponentially smaller, of order
$\frac{1}{\sqrt{R}} \e^{-2 \rho R}$; both terms are captured by the
explicit formula \eqref{eq:first-asymptotic-correction}.

Making the substitution \eqref{eq:X-specialize} in \eqref{eq:first-asymptotic-correction}
we get $R \to \infty$ asymptotics for the spectral coordinate $X_\gamma$:
\begin{equation} \label{eq:asymptotics-specialize}
X_\gamma = \exp \left[ a_\gamma R + \sum_{\mu \in \Gamma} \frac{\Omega(P_0, \mu) \IP{\gamma,\mu}}{4 \pi \I}
  \frac{\alpha_\mu + \e^{\I \vartheta}}{\alpha_\mu - \e^{\I \vartheta}} \sqrt{\frac{\pi}{R \abs{Z_\mu}}} \e^{-2 \abs{Z_\mu} R} + \delta(R, \e^{\I \vartheta}) \right],
\end{equation}
with $a_\gamma$ again given by \eqref{eq:a-gamma}.

\section{Examples} \label{sec:case-studies}

In this final section we work out the detailed statement 
of \autoref{conj:main} for two examples of polynomials $P_0(z)$.

\subsection{An \texorpdfstring{$n=2$}{n=2} example (pentagon)} \label{sec:pentagon}

We first consider a degree $2$ example,
\begin{equation} \label{eq:P0-pentagon}
  P_0(z) = \half \left(- z^2 + 1\right).
\end{equation}
The spectral curve $\Sigma$ defined in \eqref{eq:spectral-curve}
is a $3$-sheeted cover of $\C$, with ramification points of index $3$
over $z = \pm 1$.
$\Sigma$ is a one-holed torus. Thus $\Gamma$ is a lattice of rank $2$,
with nondegenerate intersection pairing.
In Figure \ref{fig:pentagon-cover-cycles} we show a convenient choice of generators 
$\gamma_1, \gamma_2$ with $\IP{\gamma_1,\gamma_2} = 1$.
\insfig{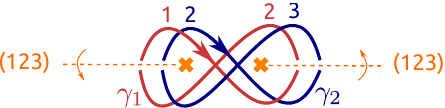}{1.6}{Generators $\gamma_1$, $\gamma_2$ for 
$\Gamma = H_1(\Sigma,\Z)$, where $P_0$ is given by \eqref{eq:P0-pentagon}. $\Sigma$
is shown as a triple cover of $\C$, with ramification 
points at $z = \pm 1$ (orange
crosses). The cover is trivialized away from branch cuts
(dashed lines). Crossing a branch
cut in the direction indicated by the arrow induces the permutation $(123)$
of sheet labels. The numerical labels next to paths indicate which
sheet of $\Sigma$ the paths lie on.}

The corresponding period integrals \eqref{eq:periods} can be computed
in closed form to give
\begin{equation}
 Z_{\gamma_1} = \e^{5 \pi \I / 6} \left( \frac{12 \times 2^{2/3} \times \pi
   ^{3/2}}{5 \Gamma
   \left(-\frac{1}{6}\right)
   \Gamma
   \left(\frac{2}{3}\right)} \right), \qquad Z_{\gamma_2} = \e^{2 \pi \I / 3} Z_{\gamma_1},
\end{equation}
or numerically
\begin{equation}
 Z_{\gamma_1} \approx -2.00324 + 1.15657 \I, \qquad Z_{\gamma_2} \approx -2.31315 \I.
\end{equation}
The analysis of spectral networks in this case (see \S\ref{sec:network-analysis-pentagon})
leads to the following:
\begin{itemize}

\item For $\abs{\vartheta} < \frac{\pi}{6}$, the spectral coordinates
attached to the network $\cW(P_0, \vartheta)$
are
\begin{equation} \label{eq:pentagon-coords}
  X_{\gamma_1} = \frac{p(1,2,3)p(3,4,5)}{p(1,3,5)p(2,3,4)}, \qquad X_{\gamma_2} = \frac{p(1,3,5)p(2,3,4)p(1,2,5)}{p(1,2,3)p(2,3,5)p(1,4,5)},
\end{equation}
where we have labeled the rays $\ell_r$ as in
Figure \ref{fig:pentagon-network} below, and we recall
the definition of Pl\"ucker coordinates
\begin{equation} \label{eq:plucker-redux}
  p(a,b,c): (v_i)_{i=1}^{n+3} \mapsto \det( v_a, v_b, v_c).
\end{equation}

\item The BPS counts are:
\begin{equation} \label{eq:bps-counts-pentagon}
  \Omega(P_0, \gamma) = \begin{cases} 1 & \text{ for } \pm \gamma \in \{\gamma_1, \gamma_2, \gamma_1 + \gamma_2 \}, \\
  0 & \text{ otherwise.} \end{cases}
\end{equation}

\end{itemize}

\insfig{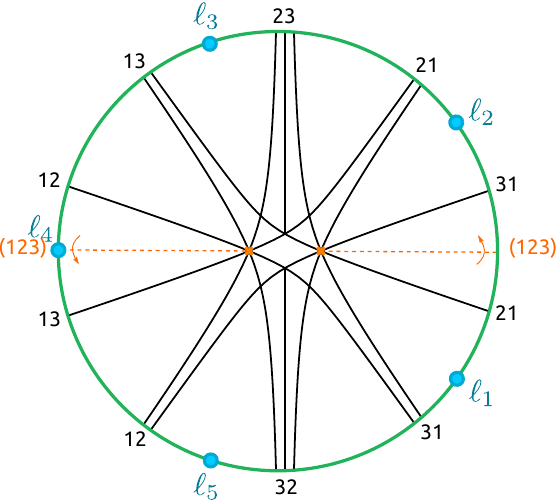}{1.35}{The spectral network $\cW(P_0, \vartheta = 0)$
where $P_0$ is given in \eqref{eq:P0-pentagon}.
It consists of $18$ WKB $\vartheta$-trajectories in all: $16$ critical trajectories emanating
from the zeroes of $P_0$, and $2$ more born from intersection points.
These $18$ trajectories approach $10$ asymptotic directions,
equally spaced around the circle at infinity.}

This is all the data necessary to formulate 
\autoref{conj:main} and its consequences 
as described in \S\ref{sec:integral-and-consequences}.
In particular, specializing the statements 
of \S\ref{sec:asymptotics} to this
example we get predictions for the $R \to \infty$
asymptotics of the $X_\gamma$. 
To give one concrete such 
prediction, we specialize
to $\vartheta = 0$ and choose $\gamma = \gamma_1$.
Then \eqref{eq:asymptotics-specialize} becomes
\begin{equation} \label{eq:asymptotics-specialize-pentagon}
  \frac{p(1,2,3)p(3,4,5)}{p(1,3,5)p(2,3,4)} = \exp \left( aR - \frac{3}{2 \sqrt{\pi \rho R}} \e^{-2 \rho R} + \delta(R) \right), \quad   a \approx -4.00648, \quad \rho \approx 2.31315,
\end{equation}
where the remainder function $\delta(R)$ obeys
\begin{equation}
 \lim_{R \to \infty} \delta(R) \sqrt{R} \e^{2 \rho R} = 0.
\end{equation}

\subsection{Spectral network analysis for the \texorpdfstring{$n=2$}{n=2} example} \label{sec:network-analysis-pentagon}

In this section we sketch the spectral network analysis
leading to \eqref{eq:pentagon-coords}
and \eqref{eq:bps-counts-pentagon}.

We begin with the spectral coordinate formulas
\eqref{eq:pentagon-coords}.
For this purpose we need to draw
the network $\cW(P_0, \vartheta = 0)$; it
is shown in Figure \ref{fig:pentagon-network}, obtained
using the Mathematica code \cite{swn-plotter}.\footnote{A version of this
code is included with the arXiv version of this paper, as {\tt swn-plotter.nb}.}
Then, according to the rules
explained in \S\ref{sec:spectral-coordinates}, we need to draw
the abelianization trees compatible with this network. 
\insfig{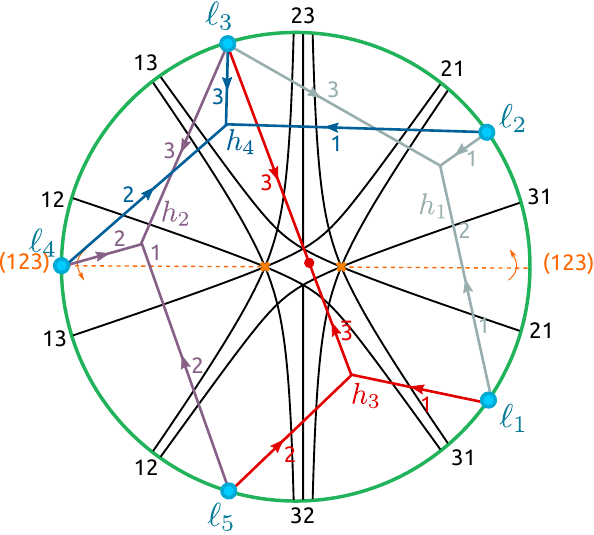}{1.4}{Four abelianization trees on the spectral
network $\cW(P_0, \vartheta = 0)$ where $P_0$ is given in \eqref{eq:P0-pentagon}.}
In Figure \ref{fig:pentagon-gamma1-trees} we show four abelianization trees $h_m$.
$h_1$, $h_2$ and $h_4$
are asymptotic abelianization trees (cf. Figure \ref{fig:asymptotic-abelianization-tree}), while $h_3$ is not asymptotic:
the legs labeled $3$ and $\overline{3}$ pass through the middle of the
spectral network. It is instructive to see why $h_3$ is compatible
with the spectral network:
the key point is that the leg labeled $3$ does not cross
any trajectory labeled $31$ or $32$ (it does cross trajectories
labeled $13$, $23$, and $21$), and the leg labeled $\overline{3}$
does not cross any trajectory labeled $13$ or $23$ (it does
cross ones labeled $32$, $31$, and $12$.)

The functions attached to these abelianization trees are
\begin{equation}
  A_{h_1} = p(1,2,3), \quad A_{h_2} = p(3,4,5), \quad A_{h_3} = p(1,3,5), \quad A_{h_4} = p(2,3,4).
\end{equation}
Combining these trees according to the rules
of \S\ref{sec:abelianization-trees}, with weights $w_1 = w_2 = +1$, $w_3 = w_4 = -1$,
we get a closed cycle which turns out to be $\gamma_1$.
Thus we have
\begin{equation}
  X_{\gamma_1} = \frac{p(1,2,3) p(3,4,5)}{p(1,3,5) p(2,3,4)}
\end{equation}
as we claimed in \eqref{eq:pentagon-coords}.
A similar construction involving $6$ abelianization trees 
gives $X_{\gamma_2}$.

In this way we obtain the formulas \eqref{eq:pentagon-coords} 
for $\vartheta = 0$. 
It still remains to see why these formulas also hold
for other $\vartheta$ with $\abs{\vartheta} < \frac{\pi}{6}$.
The reason is that the $X_\gamma$ are invariant under variations
of $\vartheta$, so long as we do not cross a phase
which is BPS-ful. (Concretely, such variations of $\vartheta$
deform $\cW(P_0, \vartheta)$ in a way which does not affect the set of
compatible abelianization trees.) 
By drawing the networks $\cW(P_0, \vartheta)$ for various
phases $\vartheta$ one can spot by eye the phases which are BPS-ful:
these are $\vartheta = \pm \frac{\pi}{6}, \pm \frac{\pi}{2}, \pm \frac{5 \pi}{6}$.\footnote{A movie containing pictures of the networks $\cW(P_0, \frac{n \pi}{300})_{n=0}^{99}$
is included with the arXiv version of this paper, as {\tt pentagon.gif}. (Recall that a shift $\vartheta \to \vartheta + \frac{\pi}{3}$ can
be compensated by a relabeling of the walls of $\cW(P_0, \vartheta)$,
so we need not explore phases beyond $[0,\frac{\pi}{3})$.)}

The BPS-ful phases are also relevant for another reason: we use
them in the process of determining the BPS counts \eqref{eq:bps-counts-pentagon}, following the rules of \S\ref{sec:bps-counts}.
At each of the BPS-ful phases a single finite web appears in
 $\cW(P_0, \vartheta)$,
consisting of a single trajectory connecting the two zeroes of $P_0$.
For example, in Figure \ref{fig:finite-web-pentagon} we show
the network $\cW(P_0, \vartheta)$ very near $\vartheta = -\frac{\pi}{6}$.
Note the two trajectories which 
almost meet head-to-head in the center of the picture.
At $\vartheta = -\frac{\pi}{6}$ these two merge into a single trajectory $p$,
whose lift $p^\Sigma$ is in homology class $-\gamma_1$.
According to the rules described in \S\ref{sec:bps-counts},
this trajectory is responsible for a nonzero count,
$\Omega(P_0, -\gamma_1) = 1$. (Note this is consistent with the fact
that $\arg Z_{-\gamma_1} = - \frac{\pi}{6}$.)
Looking similarly at the other BPS-ful phases
we get the BPS counts \eqref{eq:bps-counts-pentagon}.

\insfig{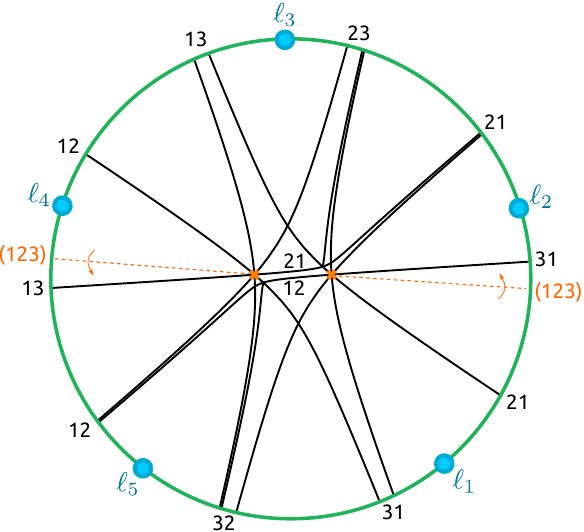}{1.25}{The spectral network $\cW(P_0, \vartheta = -\frac{\pi}{6}+0.1)$ where $P_0$ is given in \eqref{eq:P0-pentagon}. This picture can be reached by a continuous deformation of Figure \ref{fig:pentagon-network}, realized concretely
by the family of BPS-free 
networks $\cW(P_0, \vartheta)$ with $\vartheta$ varying
from $0$ to $-\frac{\pi}{6} + 0.1$.}

\subsection{An \texorpdfstring{$n=3$}{n=3} example (hexagon)} \label{sec:hexagon}

Next we consider the degree $3$ case
\begin{equation} \label{eq:P0-hexagon}
  P_0 = \frac{1}{2} \left(- z^3 + 3 z^2 + 2 \right).
\end{equation}
In this case the spectral curve $\Sigma$ defined in \eqref{eq:spectral-curve}
is a $3$-sheeted cover of $\C$, with ramification points of index $3$
over the $3$ zeroes of $P_0$.
$\Sigma$ is a $3$-holed torus. Thus $\Gamma$ is a lattice of rank $4$,
with intersection pairing of rank $2$. In Figure \ref{fig:hexagon-cover-cycles}
we show a convenient choice of generators,
with $\IP{\gamma_1,\gamma_2} = 1$ and $\gamma_3, \gamma_4$ in the kernel
of $\IP{\cdot,\cdot}$.
\insfig{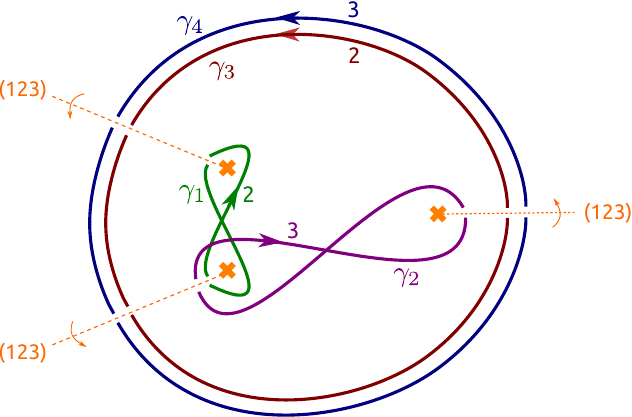}{1.2}{Generators $\{\gamma_i\}_{i=1}^4$ for $\Gamma = H_1(\Sigma,\Z)$, where $P_0$ is given by \eqref{eq:P0-hexagon}. 
The notation is as in Figure \ref{fig:pentagon-cover-cycles}.}
The numerically computed periods are
\begin{equation}
(Z_{\gamma_i})_{i=1}^4 \approx ( 2.30298,\ 5.47033+4.48792\I,\ -4.31884+2.49348\I,\ -4.98697\I ).
\end{equation}
The analysis of spectral networks in this case (see \S\ref{sec:network-analysis-hexagon}) 
leads to the following:
\begin{itemize}
\item We have
\begin{equation} \label{eq:bps-counts-hexagon}
\Omega(P_0, \gamma) = \begin{cases} 1 \text{ if } \pm \gamma \in \left\{ \begin{minipage}{7.7cm} $ (1,0,0,0),\ (0,-1,-1,-1),\ (-1,1,1,1),\\ (0,1,0,0),\ (1,-1,-1,0),\ (-1,0,1,0),\\ (0,-1,-1,0),\ (1,0,-1,-1),\ (-1,1,2,1),\\ (1,1,0,0),\ (1,-2,-2,-1),\ (-2,1,2,1) $ \end{minipage} \right\}, \\ 0 \text{ otherwise.} \end{cases}
\end{equation}
\item For $0 < \vartheta < 0.36$, the spectral coordinates
attached to the network $\cW(P_0, \vartheta)$
are
\begin{gather}
X_{\gamma_1} = \frac{q(2,3,4,5,6,1)}{p(1,5,6)p(2,3,4)}, \quad X_{\gamma_2} = \frac{p(1,5,6)p(2,3,6)p(1,4,6)}{p(1,2,6)p(1,3,6)p(4,5,6)}, \label{eq:hexagon-coords-1} \\ \quad X_{\gamma_3} = \frac{p(1,2,3)p(4,5,6)}{p(2,3,4)p(1,5,6)}, \quad X_{\gamma_4} = \frac{p(1,2,6)p(3,4,5)}{p(1,2,3)p(4,5,6)}. \label{eq:hexagon-coords-2}
\end{gather}
Here in addition to the Pl\"ucker coordinates $p(a,b,c)$ we use
the invariant (called ``hexapod invariant'' in \cite{Fomin2012})
\begin{equation} \label{eq:hexapod-invariant}
q(a,b,c,d,e,f): (v_i)_{i=1}^{n+3} \mapsto \det( v_{a} \times v_b , v_c \times v_d, v_e \times v_f).
\end{equation}
\end{itemize}

This is all the data necessary to formulate 
\autoref{conj:main} and its consequences 
as described in \S\ref{sec:integral-and-consequences}.
Unlike the case of \S\ref{sec:pentagon} above, here we get some
\ti{exact} predictions, as described in 
\S\ref{sec:exact}: indeed, since $\gamma_3$ and $\gamma_4$
are in the kernel of $\IP{\cdot,\cdot}$ we get exact
formulas for $X_{\gamma_3}$ and $X_{\gamma_4}$.
For example, suppose we specialize
to $\vartheta = 0.2$ and choose $\gamma = \gamma_3$.
Then \eqref{eq:exact-specialized} becomes
\begin{equation} \label{eq:exact-specialize-hexagon}
 \frac{p(1,2,3)p(4,5,6)}{p(2,3,4)p(1,5,6)} = \exp \left( aR \right), \quad   a \approx -7.4748. 
 \end{equation}
We also get asymptotic predictions
for the other coordinates, as we did in \S\ref{sec:pentagon} above.
For example, specializing to $\vartheta = 0.2$ and $\gamma = \gamma_1$,
\eqref{eq:asymptotics-specialize} becomes
\begin{equation}
  \frac{q(2,3,4,5,6,1)}{p(1,5,6)p(2,3,4)} = \exp \left( aR + \frac{c}{\sqrt{R}} \e^{-2 \rho R} + \delta(R) \right)
\end{equation}
where
\begin{equation}
  a \approx 4.5142, \quad \rho \approx 2.3030, \quad c \approx 0.1961,
\end{equation}
and the remainder function $\delta(R)$ obeys
\begin{equation}
 \lim_{R \to \infty} \delta(R) \sqrt{R} \e^{2 \rho R} = 0.
\end{equation}

\subsection{Spectral network analysis for the \texorpdfstring{$n=3$}{n=3} example} \label{sec:network-analysis-hexagon}

In this section we sketch the spectral network analysis leading
to \eqref{eq:bps-counts-hexagon}, \eqref{eq:hexagon-coords-1}, \eqref{eq:hexagon-coords-2}.
This is parallel
to \S\ref{sec:network-analysis-pentagon} above.

The spectral coordinates \eqref{eq:hexagon-coords-1},
\eqref{eq:hexagon-coords-2}
are obtained as follows.
We first draw the network $\cW(P_0, \vartheta = 0.1)$, 
shown in Figure \ref{fig:hexagon-network}, obtained
using \cite{swn-plotter}.
\insfig{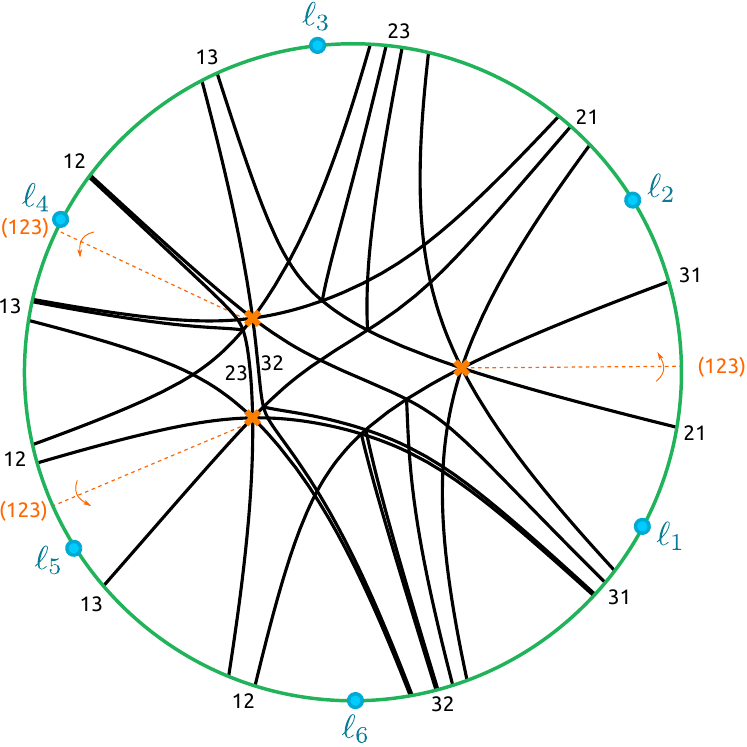}{1.20}{The spectral network $\cW(P_0, \vartheta = 0.1)$ where $P_0$ is given in \eqref{eq:P0-hexagon}. It consists of $31$ WKB $\vartheta$-trajectories: $24$ critical trajectories emanating from the zeroes of $P_0$, and $7$ more born from intersection points.}
\insfig{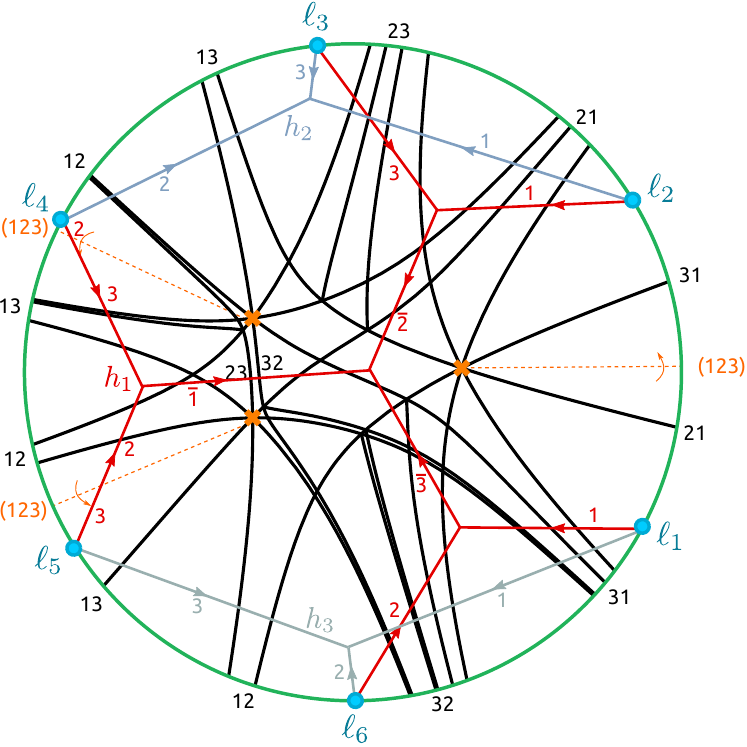}{1.20}{Three abelianization 
trees compatible with the 
spectral network $\cW(P_0, \vartheta = 0.1)$ where $P_0$ is given
in \eqref{eq:P0-hexagon}.}
Then, we draw compatible abelianization trees.
In Figure \ref{fig:hexagon-gamma1-trees} we show three abelianization trees $h_m$.
$h_2$ and $h_3$
are asymptotic abelianization trees (cf. Figure \ref{fig:asymptotic-abelianization-tree})
while $h_1$ is not. $h_1$ is a more interesting tree than we have seen
thus far: it is a ``hexapod'' in the terminology of \cite{Fomin2012}.
The functions attached to these abelianization trees are
\begin{equation}
  A_{h_1} = q(2,3,4,5,6,1), \qquad A_{h_2} = p(1,5,6), \qquad A_{h_3} = p(2,3,4).
\end{equation}
Combining these according to the rules
of \S\ref{sec:abelianization-trees}, with weights $w_1 = +1$, $w_2 = w_3 = -1$,
we get a closed cycle which turns out to be $\gamma_1$.
Thus we have
\begin{equation}
  X_{\gamma_1} = \frac{q(2,3,4,5,6,1)}{p(1,5,6) p(2,3,4)}
\end{equation}
as we claimed in \eqref{eq:hexagon-coords-1}. Similar constructions
and computations give the other $X_{\gamma_i}$ of
\eqref{eq:hexagon-coords-1}, \eqref{eq:hexagon-coords-2}.

The BPS-ful phases and the BPS counts $\Omega(P_0, \gamma)$ in \eqref{eq:bps-counts-hexagon}
were determined, as in the $n=2$ case above, by direct 
examination of the networks
$\cW(P_0, \vartheta)$ as $\vartheta$ varies.\footnote{A movie containing pictures of the networks $\cW(P_0, \frac{n \pi}{300})_{n=0}^{99}$
is included with the arXiv version of this paper, as {\tt hexagon.gif}.}
For each unordered pair $(z,z')$ of zeroes of $P_0$, we find $6$ finite BPS webs
which are single trajectories connecting $z$ to $z'$. For example, at
$\vartheta = 0$ we get a finite BPS web of charge $\gamma_1$, 
connecting the two leftmost
zeroes; indeed, in Figure \ref{fig:hexagon-network} one can easily see
two trajectories in $\cW(P_0, \vartheta = 0.1)$
which almost meet head-on, and which do
meet head-on when $\vartheta$ is adjusted from $0.1$ to $0$.
Similarly, at $\vartheta \approx 0.36$ we get a finite BPS
web of charge $\gamma_1 - \gamma_3 - \gamma_4$, consisting of
a single trajectory connecting the zeroes at top left and at
right. There are also $6$
more finite BPS webs which are three-string junctions involving all three 
zeroes of $P_0$; these webs have the topology shown in the 
right side of Figure \ref{fig:finite-webs}. For each of the $24$
finite webs we worked out the corresponding charge $\gamma$ following
the rules of \S\ref{sec:bps-counts}.
This gives
the $24$ nonzero $\Omega(P_0, \gamma)$ listed in
\eqref{eq:bps-counts-hexagon}.

\bibliography{harmonic-predictions-paper}

\end{document}